\documentclass[11pt]{article}
\usepackage{fullpage}
\usepackage{authblk}
\usepackage{natbib}
\usepackage{algorithm}
\usepackage{algpseudocode}
\usepackage{amsmath,amsthm,amsfonts}
\usepackage{amsmath}
\usepackage[subnum]{cases}
\usepackage{tcolorbox}
\usepackage{mathrsfs}
\usepackage{amssymb}
\usepackage{color}
\usepackage{mathrsfs}
\usepackage{empheq}
\usepackage{enumitem}
\usepackage{bm}
\usepackage{multirow}
\usepackage{booktabs}
\usepackage{makecell}
\usepackage{graphicx}
\usepackage{subcaption}
\usepackage{comment}
\usepackage{cases}
\usepackage{appendix}
\usepackage{tikz}
\usetikzlibrary{arrows,shapes}
\usepackage[colorlinks= true, linkcolor=red, citecolor=blue, urlcolor=black]{hyperref}
\usepackage{cleveref}
\usepackage{marginnote}
\usepackage{diagbox} 
\numberwithin{equation}{section}

\theoremstyle{definition}

\newtheorem{theorem}{Theorem}[section]
\newtheorem{lemma}[theorem]{Lemma}
\newtheorem{corollary}[theorem]{Corollary}
\newtheorem{remark}[theorem]{Remark}

\newtheorem{defn}[theorem]{Definition}

\newcommand{\mr}{\mathbb{R}}

\usepackage{booktabs}
\usepackage{tikz}
\usepackage{pgfplots}
\usetikzlibrary{calc,intersections}

\usetikzlibrary{shapes.geometric, arrows, positioning}

\tikzset{
	mybox/.style  = {draw, rectangle, minimum width=4cm, minimum height=0.8cm, text centered, text width=4.4cm,   
		font=\normalsize},
	box/.style  = {draw, rectangle, minimum width=2.0cm, minimum height=0.6cm, text centered, text width=3.0cm,   
		font=\normalsize},
	myarrow/.style = {line width=0.2pt, draw=black, -triangle 60, postaction={draw, line width=0.2pt, shorten >=10pt,-}}
}

\tikzstyle{arrow} = [->, >=stealth, -triangle 60]

\allowdisplaybreaks

\makeatletter
\newcommand{\leqnomode}{\tagsleft@true}
\newcommand{\reqnomode}{\tagsleft@false}
\makeatother

\begin{document}
%\maketitle

\title{A Family of Controllable Momentum Coefficients for Forward-Backward Accelerated Algorithms}
\author[1]{Mingwei Fu}
\author[2,3]{Bin Shi\thanks{Corresponding author: \url{binshi@fudan.edu.cn} } }
\affil[1]{School of Mathematical Sciences, University of Chinese Academy of Sciences, Beijing 100049, China}
\affil[2]{Shanghai Institute for Mathematics and Interdisciplinary Sciences, Shanghai 200433, China}
\affil[3]{Center for Mathematics and Interdisciplinary Sciences, Fudan University, Shanghai 200433, China}
\date\today

\maketitle

%\bin{Vanishing Viscosity should be included in the main part, the view is key important. At least, it is mentioned a little.}\wjs{OK.}
\begin{abstract}
Nesterov’s accelerated gradient method (\texttt{NAG}) marks a pivotal advancement in gradient-based optimization, achieving faster convergence compared to the vanilla gradient descent method for convex functions. However, its algorithmic complexity when applied to strongly convex functions remains unknown, as noted in the comprehensive review by~\citet{chambolle2016introduction}. This issue, aside from the critical step size, was addressed by~\citet{li2024linear2}, with the monotonic case further explored by~\citet{fu2024lyapunov}. In this paper, we introduce a family of controllable momentum coefficients for forward-backward accelerated methods, focusing on the critical step size $s=1/L$. Unlike traditional linear forms, the proposed momentum coefficients follow an $\alpha$-th power structure, where the parameter $r$ is adaptively tuned to $\alpha$. Using a Lyapunov function specifically designed for $\alpha$, we establish a controllable $O\left(1/k^{2\alpha} \right)$ convergence rate for the \texttt{NAG}-$\alpha$ method, provided that $r > 2\alpha$. At the critical step size, \texttt{NAG}-$\alpha$ achieves an inverse polynomial convergence rate of arbitrary degree by adjusting $r$ according to $\alpha > 0$. We further simplify the Lyapunov function by expressing it in terms of the iterative sequences $x_k$ and $y_k$, eliminating the need for phase-space representations. This simplification enables us to extend the controllable $O \left(1/k^{2\alpha} \right)$ rate to the monotonic variant, \texttt{M-NAG}-$\alpha$, thereby enhancing optimization efficiency. Finally, by leveraging the fundamental inequality for composite functions, we extended the controllable $O\left(1/k^{2\alpha} \right)$ rate to proximal algorithms, including the fast iterative shrinkage-thresholding algorithm (\texttt{FISTA}-$\alpha$) and its monotonic counterpart (\texttt{M-FISTA}-$\alpha$).
\end{abstract}

%\thispagestyle{empty}
%\setcounter{page}{0}
%
%%\new-page
%\tableofcontents
%%\newpage

\section{Introduction}
\label{sec: intro}

Over the past few decades, machine learning has experienced unprecedented growth, driving transformative advancements across a wide range of scientific disciplines and industrial applications. At the heart of this revolution lies the fundamental challenge of unconstrained optimization, particularly the task of minimizing a convex objective function without constraints. Mathematically, this problem can be expressed as:
\[
\min_{x \in \mathbb{R}^d} f(x).
\]
Addressing this problem has been a central focus of optimization research, leading to the development of gradient-based algorithms, which form the backbone of modern machine learning methodologies. These algorithms are celebrated for their computational efficiency and low memory requirements, making them particularly well-suited for large-scale problems, where computational and memory resources are often limited. Beyond their practical utility, gradient-based methods have also catalyzed profound theoretical advancements while enabled impactful real-world applications. As a result, these algorithms have now become indispensable tools in both academic research and industrial practice.

Among the family of gradient-based methods, the vanilla gradient descent method stands out as one of the simplest yet most fundamental approaches, which can be traced back to to the pioneering work of~\citet{cauchy1847methode}. Starting from an initial point $x_0 \in \mathbb{R}^d$, it follows a straightforward iterative update rule: 
\[
x_{k} = x_{k-1} - s \nabla f(x_{k-1}),
\]
where $s>0$ denotes the step size. For convex functions, the vanilla gradient descent method guarantees convergence with the following rate:
\[
f(x_k) - f(x^{\star}) \leq O\left( \frac{1}{sk} \right), 
\]
where $x^{\star}$ denotes the global minimum point of $f(x)$. In the mid-1980s,~\citet{nesterov1983method} introduced a groundbreaking innovation that revolutionized optimization theory: the accelerated gradient method, now widely known as Nesterov’s accelerated gradient method (\texttt{NAG}). This method incorporates an elegant two-step mechanism to achieve faster convergence. Given an initial point $x_0 = y_0 \in \mathbb{R}^d$,~\texttt{NAG} is implemented through the following iterative update rules:
\[
\left\{ \begin{aligned}
         & x_{k} = y_{k-1} - s \nabla f(y_{k-1}), \\
         & y_{k} = x_{k} + \frac{k-1}{k+r} \left( x_{k} - x_{k-1} \right), 
         \end{aligned} \right.
\]
where $s>0$ is the step size, and $r \geq 2$ is the momentum coefficient. Compared to the vanilla gradient descent,~\texttt{NAG} achieves an accelerated convergence rate:
\[
f(x_k) - f(x^{\star}) \leq O\left( \frac{1}{sk^2} \right). 
\]
This acceleration represents a major leap in optimization efficiency, particularly for large-scale problems. To extend the applicability to composite optimization problems, especially those arising in image science,~\citet{beck2009fast} introduced a proximal variant of~\texttt{NAG}, known as the fast iterative shrinkage-thresholding algorithm (\texttt{FISTA}). This approach was later refined by~\citet{beck2017first}, who introduced monotonic versions of both~\texttt{NAG} and~\texttt{FISTA}.  Collectively, these algorithms are referred to as forward-backward accelerated gradient algorithms, highlighting their forward-backward structure. 

Despite their solid theoretical foundations, these forward-backward accelerated gradient algorithms have been primarily designed for convex functions. However, many real-world applications often involve strongly convex but ill-conditioned, characterized by Hessian eigenvalues bounded away from zero. For such strongly convex functions, the forward-backward accelerated gradient algorithms have yet to demonstrate superior performance over the vanilla gradient descent and its variants. Indeed, a critical question remains unsolved: whether these forward-backward algorithms can achieve linear convergence for strongly convex functions. This challenge was underscored in the seminal review by~\citet{chambolle2016introduction}. Recent efforts have sought to address this gap. Building on the high-resolution differential equation framework~\citep{shi2022understanding},~\citet{li2024linear} analyzed forward-backward accelerated method for step sizes $0 < s < 1/L$, demonstrating linear convergence.\footnote{The parameter $L$ is the Lipschitz constant, which is defined in detail in~\Cref{sec: prelim}.} Additionally,~\citet{fu2024lyapunov} developed a novel Lyapunov function, which extends linear convergence to monotonic variants of these algorithms. Nevertheless, a critical limitation persists: for the critical step size $s=1/L$, no theoretical guarantees currently exist for the forward-backward accelerated gradient algorithms. This unresolved issue raises two fundamental questions in the context of strongly convex functions and the critical step size $s = 1/L$:
\begin{tcolorbox}
\begin{itemize}
\item How fast can forward-backward accelerated gradient algorithms converge? 
\item Can the forward-backward accelerated gradient algorithms be improved to achieve faster convergence? 
\end{itemize}
\end{tcolorbox}
\noindent Answering these questions represents a pivotal challenge in optimization theory, with the potential to further bridge the gap between theoretical developments and practical applications in machine learning.

%Answering these questions remains a key challenge in optimization theory and has the potential to further bridge the gap between theory and practice in machine learning.

%%%%%%%%%%%%%%%%%%%%%%%%%%%%%%%%%%%%%%%%%%%%%%%%%%%%%%%%%%%%%%%%%%%%%%%%%%%%%%%%%%%%%%%%%%%%
\subsection{Motivation: a family of controllable momentum coefficients}
\label{subsec: motivation}

The classical~\texttt{NAG} method leverages a momentum coefficient, which is typically expressed as: 
\[
\frac{k-1}{k+r},
\]  
where $k$ denotes the iteration number, and $r$ is a constant that governs the strength of the momentum term. This simple yet effective formulation has consistently demonstrated reliable performance across various convex optimization problems. However, its inherent simplicity imposes a notable limitation: the rigidity of this momentum coefficient restricts the algorithm's ability to adapt to problem-specific structures. This limitation becomes particularly pronounced in scenarios involving strongly convex or ill-conditioned objective functions, where a more flexible approach could better exploit the underlying structure for improved performance.

To address these limitations, we introduce a more generalized and flexible momentum coefficient. Specifically, we extend the classical linear form to an $\alpha$-th power formulation: 
\begin{equation}
\label{eqn: alpha-moemtnum-coefficient}
\frac{(k-1)^{\alpha}}{k^{\alpha} + rk^{\alpha-1}},
\end{equation}
where $\alpha$ is a tunable parameter that controls the degree of acceleration, and $r$ is a controllable parameter that can be adjusted in accordance with the value of $\alpha$. This generalized momentum coefficient introduces enhanced adaptability, enabling the algorithm to dynamically tailor its momentum dynamics to the specific characteristics of the problem. By leveraging this flexibility, the proposed formulation has the potential to achieve faster convergence rates, particularly in challenging optimization landscapes.

%This generalized momentum coefficient provides enhanced adaptability, enabling the algorithm to potentially achieve faster convergence rates by tailoring the momentum dynamics to specific problem characteristics.

%This generalized momentum term significantly enhances the algorithm's flexibility, allowing it to better adopt to the specific requirements of different optimization problems.

Building upon this new momentum coefficient~\eqref{eqn: alpha-moemtnum-coefficient}, we propose an extended variant of the classical~\texttt{NAG} method, termed~\texttt{NAG}-$\alpha$. The forward-backward update rules for this generalized algorithm, incorporating the new momentum coefficient from~\eqref{eqn: alpha-moemtnum-coefficient}, are as follows: 
\begin{subequations}
\label{eqn: nag-alpha}
    \begin{empheq}[left=\empheqlbrace]{align}
         & x_{k} = y_{k-1} - s \nabla f(y_{k-1})                                                                                           \label{eqn: nag-alpha-gradient}         \\
         & y_{k} = x_{k} + \frac{(k-1)^{\alpha}}{k^{\alpha}+rk^{\alpha-1}} \left( x_{k} - x_{k-1} \right)         \label{eqn: nag-alpha-momentum} 
    \end{empheq}
\end{subequations}
where $s>0$ is the step size. This extension allows for greater control over the acceleration dynamics through the parameter $\alpha$, making the algorithm better suited for a broader range of optimization problems. Inspired by prior works~\citep{beck2017first, fu2024lyapunov}, we further extend the momentum coefficient~\eqref{eqn: alpha-moemtnum-coefficient} to construct a monotonic variant of~\texttt{NAG}-$\alpha$, which we denote as \texttt{M-NAG}-$\alpha$. This variant incorporates an adaptive mechanism to ensure monotonicity in the progression of the function values, enhancing the robustness of the method.  The update rules for \texttt{M-NAG}-$\alpha$ are defined as follows:
\begin{subequations}
\begin{empheq}[left=\empheqlbrace]{align}
& z_{k-1} = y_{k-1} - s\nabla f(y_{k-1}),                                                                                                                \label{eqn: m-nag-gradient}       \\
& x_{k} = \left\{ \begin{aligned} 
                          & z_{k-1}, \quad \text{if}\; f(z_{k-1}) \leq f(x_{k-1}),  \\
                          & x_{k-1}, \quad \text{if}\; f(z_{k-1}) > f(x_{k-1}),
                          \end{aligned} \right.                                                                                                                    \label{eqn: m-nag-comparison}   \\    
& y_{k} = x_{k} + \frac{(k-1)^{\alpha}}{k^{\alpha}+rk^{\alpha-1}} \left( x_{k} - x_{k-1} \right) + \frac{(k-1)^{\alpha} + r(k-1)^{\alpha-1}}{k^{\alpha}+rk^{\alpha-1}} \left(z_{k-1} - x_{k} \right),            \label{eqn: m-nag-momentum}
\end{empheq}
\end{subequations}
where $s>0$ is the step size. This monotonic extension integrates momentum-based acceleration with function-value monotonicity, thereby improving the algorithm's versatility and ensuring more stable convergence behavior.

%This flexibility opens the door to faster convergence and improved performance across a broader range of optimization tasks, especially for challenging problems with strong convexity and ill-conditioning.

\begin{figure}[htb!]
\centering
\begin{subfigure}[t]{0.475\linewidth}
\centering
\includegraphics[scale=0.16]{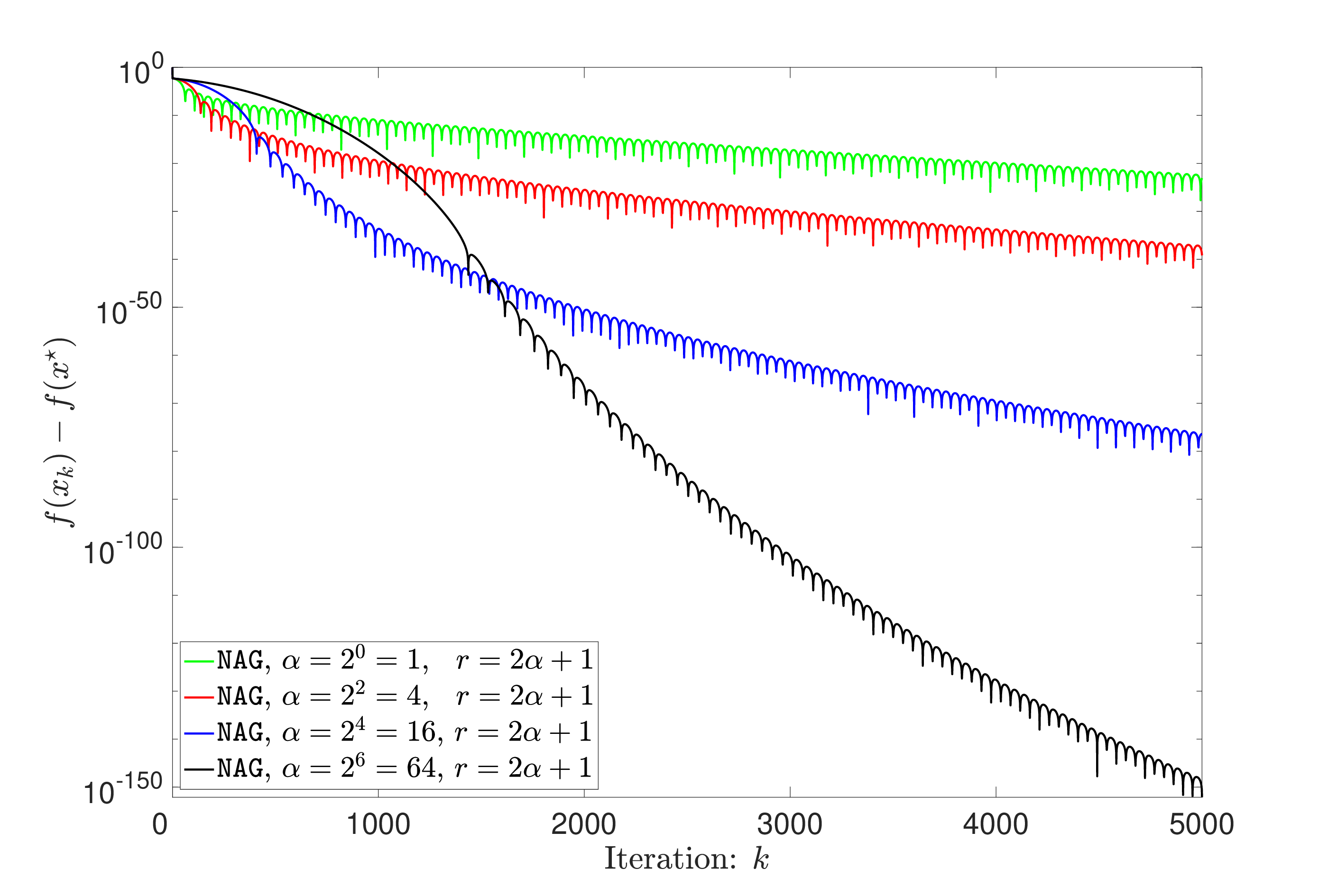}
\caption{\texttt{NAG}-$\alpha$}
\label{subfig: nag}
\end{subfigure}
\begin{subfigure}[t]{0.475\linewidth}
\centering
\includegraphics[scale=0.16]{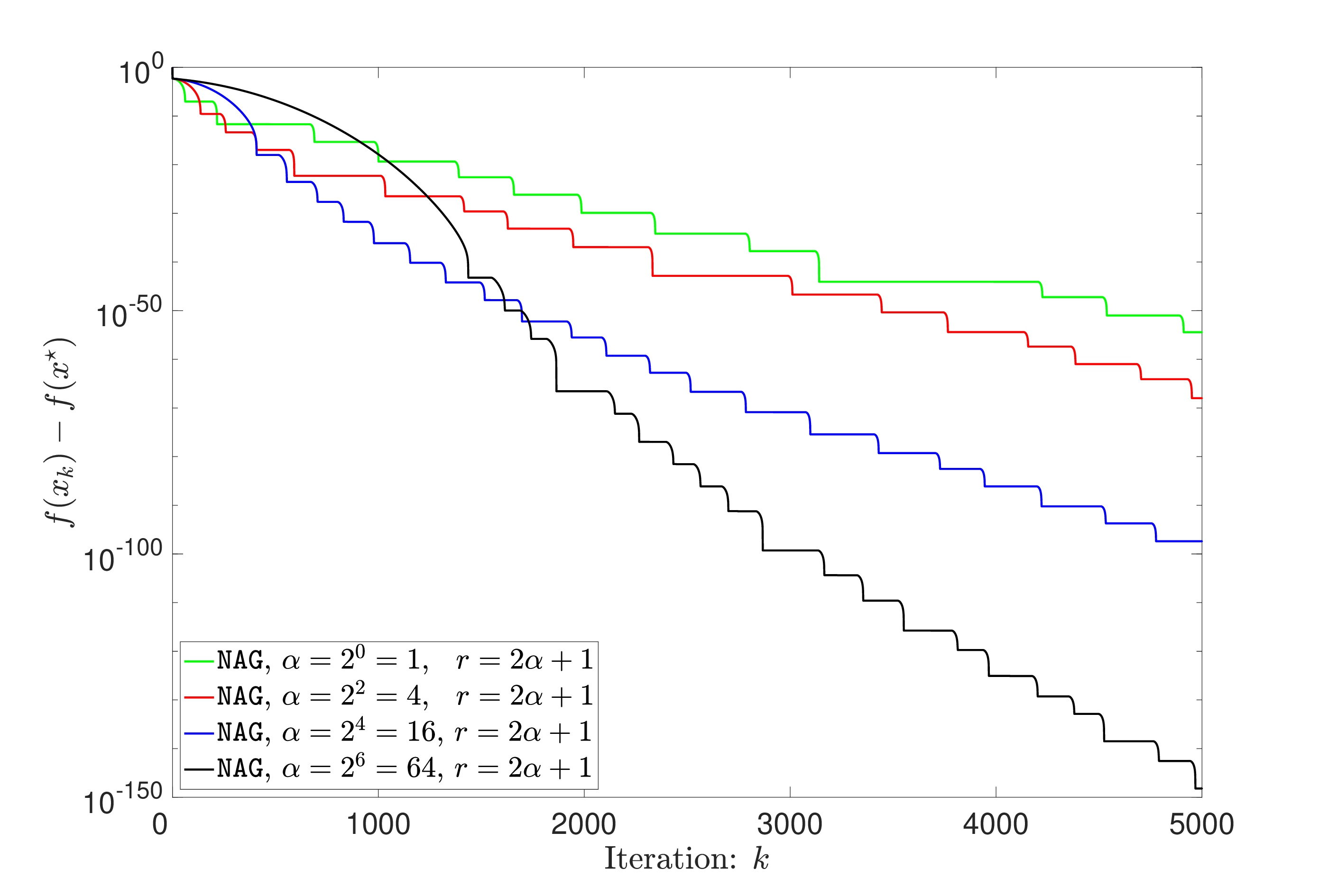}
\caption{\texttt{M-NAG}-$\alpha$}
\label{subfig: mnag}
\end{subfigure}
%\begin{subfigure}[t]{0.325\linewidth}
%\centering
%\includegraphics[scale=0.3]{fig/recovery.eps}
%\caption{Deblurred and Denoisy}
%\end{subfigure}
%\caption{Deblurring and Denoising an image of an elephant by ISTA and FISTA.} 
\caption{Iterative progression of function values for~\texttt{NAG}-$\alpha$ and~\texttt{M-NAG}-$\alpha$ applied to the quadratic function $f(x_1, x_2) = 5 \times 10^{-3}x_1^2 + x_2^2$. The experiments are performed with $s = 1/L = 0.5$.}
\label{fig: nag-alpha}
\end{figure}

To evaluate the practical performance of~\texttt{NAG}-$\alpha$ and~\texttt{M-NAG}-$\alpha$, we conduct numerical experiments on a quadratic function under the critical step size.~\Cref{fig: nag-alpha} illustrates the iterative progression of function values for both algorithms, highlighting the impact of varying $\alpha$.
The numerical results demonstrate that, as $\alpha$ increases and $r=2\alpha+1$ is appropriately adjusted, both~\texttt{NAG}-$\alpha$ and~\texttt{M-NAG}-$\alpha$ exhibit progressively faster convergence rates. These observations underscore the potential of the generalized momentum coefficient to significantly enhance algorithmic efficiency, leading to the following key conjectures:
\begin{tcolorbox}
\begin{itemize}
\item For strongly convex functions, the convergence rates of both~\texttt{NAG}-$\alpha$ and~\texttt{M-NAG}-$\alpha$ vary with the parameter $\alpha$, provided that the control parameter $r$ is appropriately tuned. 
\item At the critical step size $s=1/L$, both~\texttt{NAG}-$\alpha$ and~\texttt{M-NAG}-$\alpha$ are expected to offer theoretical guarantees for faster convergence compared to the classical~\texttt{NAG} method and its monotonic variant. 
\end{itemize}
\end{tcolorbox}
%\noindent This novel formulation and its extensions open promising avenues for future research, particularly in developing optimization methods that dynamically adapt to problem-specific characteristics, thereby improving their overall robustness and versatility.
%\noindent  The flexibility introduced by the $\alpha$-dependent momentum term allows \texttt{NAG}-$\alpha$  to dynamically adjust to varying levels of acceleration. This adaptability enhances its versatility, making it a powerful and highly flexible method capable of efficiently tackling a wide range of optimization problems with diverse convergence behaviors.

%This adaptability in the momentum term allows \texttt{NAG}-$\alpha$ to dynamically adjust to different levels of acceleration, thereby making it a highly adaptable and versatile method. As a result, it is well-suited to address a broad range of optimization problems with varying convergence characteristics.

%This adaptability allows the algorithm to effectively tune itself to different levels of acceleration, highlighting its potential for a wide range of optimization problems. It is observed that, with the value of the parameter $\alpha$ increasing and the controllable parameter $r=2\alpha+1$,~\texttt{NAG}-$\alpha$ converges faster and faster at the critical step size $s=1/L$. Hence, we can conjecture that 

% The parameter $\alpha$ takes varying value, and the controllable parameter is set to $r=2\alpha+1$. 

%%%%%%%%%%%%%%%%%%%%%%%%%%%%%%%%%%%%%%%%%%%%%%%%%%%%%%%%%%%%%%%%%%%%%%%%%%%%%%%%%%%%%%%%%%%%
\subsection{Overview of contributions}
\label{subsec: contributions}

In this paper, we introduce a novel family of momentum coefficients for forward-backward accelerated methods and thoroughly investigate their convergence behavior, with a particular focus on the critical step size $s = 1/L$. Instead of the conventional linear form, the momentum coefficient takes the $\alpha$-th power form, where the parameter $r$ is both controlled and adaptively tuned to the value of $\alpha$. 

\begin{itemize}
\item[(\textbf{I})]  We begin by constructing a novel Lyapunov function tailored to the momentum parameter $\alpha$, which intentionally excludes the kinetic energy term. This formulation allows us to derive a controllable $O\left( 1/k^{2\alpha} \right)$ convergence rate for the \texttt{NAG}-$\alpha$ method, provided that the parameter $r$ satisfies the condition $r > 2\alpha$. At the critical step size $s=1/L$, the~\texttt{NAG}-$\alpha$ method achieves an inverse polynomial convergence rate with arbitrary degree by adaptively tuning the parameter $r$ in accordance with the parameter $\alpha > 0$. 

\item[(\textbf{II})]  The Lyapunov function that we construct for~\texttt{NAG}-$\alpha$, by omitting the kinetic energy term, provides a simplified way to express the kinetic energy directly in terms of the iterative sequences $x_k$ and $y_k$, avoiding the need for phase-space representations. This key observation enables us to extend the controllable $O\left( 1/k^{2\alpha} \right)$ convergence rate from \texttt{NAG}-$\alpha$ to its monotonic variant, \texttt{M-NAG}-$\alpha$, offering an efficient method for accelerated optimization.

\item[(\textbf{III})]  Building on the fundamental inequality for composite functions derived in~\citep{li2024linear2}, which serves as the proximal counterpart to the strongly convex inequality, we extend the controllable $O\left( 1/k^{2\alpha}\right)$ convergence rate to the two proximal algorithms, the fast iterative shrinkage thresholding algorithm (\texttt{FISTA}-$\alpha$) and its monotonic counterpart (\texttt{M-FISTA}-$\alpha$).
\end{itemize}

%%%%%%%%%%%%%%%%%%%%%%%%%%%%%%%%%%%%%%%%%%%%%%%%%%%%%%%%%%%%%%%%%%%%%%%%%%%%%%%%%%%%%%%%%%%%
\subsection{Related works and organization}
\label{subsec: related}

The introduction of~\texttt{NAG}, originally proposed by~\citet{nesterov1983method}, represents a landmark moment in the evolution of gradient-based optimization algorithms.  This two-step forward-backward algorithm fundamentally transformed optimization by incorporating acceleration into the traditional gradient descent method, resulting in significant improvements in convergence rates. Building upon this breakthrough,~\citet{beck2009fast} introduced a proximal version of the fundamental inequality, which paved the way for the development of~\texttt{FISTA}. This new algorithm extended the power of \texttt{NAG}  to handle composite objective functions,  making it a more efficient tool for a wide range of applications, including signal processing and image reconstruction.~\citet{beck2017first} further advanced this field by proposing a monotonically convergent version of the forward-backward accelerated gradient algorithms. This variant ensured that the performance of the algorithm would consistently improve over time. Simultaneously, adaptive restart strategies were explored to maintain monotonicity and enhance the robustness of forward-backward accelerated algorithms, with some contributions from the works of~\citet{giselsson2014monotonicity} and~\citet{o2015adaptive} .

Over the past decade, there has been growing interest in understanding the acceleration phenomenon in optimization algorithms through the lens of differential equations. This exploration began with works~\citep{attouch2012second, attouch2014dynamical}, which sparked a surge of research into the dynamics underlying accelerated methods. The topic gained further prominence with the development of a low-resolution differential equation for modeling~\texttt{NAG}~\citep{su2016differential}, a variational perspective on acceleration~\citep{wibisono2016variational}, and studies focusing on the faster convergence rate of function values~\citep{attouch2016rate}.  The acceleration mechanism was ultimately clarified through comparisons between~\texttt{NAG} method with Polyak's heavy-ball method, with pivotal insights emerging from the high-resolution differential equation framework introduced by~\citet{shi2022understanding}. This framework revealed that gradient correction is effectively achieved through an implicit velocity update, providing a more comprehensive understanding of how acceleration arises. Subsequent works~\citep{chen2022gradient, chen2022revisiting} further demonstrated that this framework is particularly well-suited for \texttt{NAG}. Buildling on this insights, significant progress was made by extending the framework to address composite optimization problems, encompassing both convex and strongly convex functions~\citep{li2022proximal, li2022linear}. This extension was achieved by refining the proximal inequality and generalizing the framework to accommodate a wider range of optimization scenarios, including the overdamped case, as shown in~\citep{chen2023underdamped}. These advancements have greatly enhanced our understanding of acceleration mechanisms and paved the way for more efficient, robust, and versatile optimization algorithms.

The remainder of this paper is structured as follows.~\Cref{sec: prelim} introduces the definitions and fundamental inequalities for strongly convex functions and composite functions, laying the groundwork for the analysis.~\Cref{sec: nag} outlines the construction of the Lyapunov function and establishes the controllable $O(1/k^{2\alpha})$ convergence rate for both~\texttt{NAG}-$\alpha$ and~\texttt{FISTA}-$\alpha$.~\Cref{sec: m-nag-fista} extends he controllable $O(1/k^{2\alpha})$ convergence rate  to the monotonically accelerated forward-backward algorithms,~\texttt{M-NAG}-$\alpha$ and \texttt{M-FISTA}-$\alpha$, leveraging the novel Lyapnov function that excludes the kinetic energy. Finally,~\Cref{sec: conclusion} concludes this papers and proposes potential avenues for future research.

\section{Preliminaries}
\label{sec: prelim} 

In this paper, we adopt notations that align closely with the conventions established in~\citep{nesterov2018lectures, shi2022understanding, li2024linear2}, while incorporating slight modifications tailored to the specific requirements of our analysis. Let $\mathcal{F}^0(\mathbb{R}^d)$ denote the class of continuous convex functions defined on $\mathbb{R}^d$. Specifically, a continuous function $g \in \mathcal{F}^{0}(\mathbb{R}^d)$ if it satisfies the convex condition:
\[
g\left( \alpha x + (1 - \alpha)y \right) \leq \alpha g(x) + (1 - \alpha)g(y), \quad \forall x,y \in \mathbb{R}^d,\; \alpha \in [0, 1].
\] 
Within this class $\mathcal{F}^0(\mathbb{R}^d)$, we define the subclass $\mathcal{F}^1_L(\mathbb{R}^d)$, which consists of continuously differentiable functions with Lipschitz continuous gradients. Specifically, a function $f \in \mathcal{F}^{1}_{L}(\mathbb{R}^d)$ if $f \in \mathcal{F}^{0}(\mathbb{R}^d)$ is differentiable, and satisfies the Lipschitz gradient condition:
\begin{equation}
\label{eqn: grad-lip}
\| \nabla f(x) - \nabla f(y) \| \leq L \| x - y \|, \quad \forall x,y \in \mathbb{R}^d.
\end{equation}
We further narrow our focus to the class $\mathcal{S}_{\mu,L}^{1}(\mathbb{R}^d)$, a subclass of $\mathcal{F}^1_L(\mathbb{R}^d)$, where each function is $\mu$-strongly convex for some $0 < \mu \leq L$. Specifically, a function $f \in \mathcal{S}_{\mu,L}^{1}(\mathbb{R}^d)$ if $f \in \mathcal{F}^1_L(\mathbb{R}^d)$ and satisfies the $\mu$-strongly convex condition:
\begin{equation}
\label{eqn: defn-scvx}
f(y) \geq f(x) +  \langle \nabla f(y), y - x \rangle + \frac{\mu}{2} \|y - x\|^2, \quad \forall x,y \in \mathbb{R}^d.
\end{equation}
For any $f \in \mathcal{S}_{\mu,L}^{1}(\mathbb{R}^d)$, the following fundamental inequality holds for any step size  $s \in (0, 1/L]$:
\begin{equation}
\label{eqn: fund-inq-smooth}
f(y - s\nabla f(y)) - f(x) \leq \left\langle \nabla f(y), y - x \right\rangle - \underbrace{\frac{\mu}{2} \|y - x\|^2}_{\mathbf{SC}} - \frac{s}{2} \|\nabla f(y)\|^2, \quad \forall x,y \in \mathbb{R}^d.
\end{equation}
%for any $x,y \in \mathbb{R}^d$. 

%Additionally, we denote $x^\star$ as the unique minimizer of the objective function. 

%Furthermore, for any $f \in \mathcal{S}_{\mu,L}^{1}(\mathbb{R}^d)$, we have the following inequality 
%\begin{equation}
%\label{eqn: key-inq-smooth}
%\left\| \nabla f(y) \right\|^2 \geq 2\mu \left( f(y - s\nabla f(y)) - f(x^{\star}) \right),
%\end{equation}
%for any $y \in \mathbb{R}^d$. 

We now turn our attention to composite functions of the form $\Phi= f + g$, where $f \in \mathcal{S}_{\mu,L}^1(\mathbb{R}^d)$ and $g \in \mathcal{F}^0(\mathbb{R}^d)$. Building upon the methodologies in~\citep{beck2009fast, su2016differential}, we introduce the concept of the $s$-proximal value, which serves as a cornerstone of our subsequent analysis.
\begin{defn}[$s$-Proximal Value]
\label{defn: proximal-value}
Let the step size satisfy $s \in (0, 1/L]$. For any $f\in \mathcal{S}_{\mu,L}^1(\mathbb{R}^d)$ and $g \in \mathcal{F}^0(\mathbb{R}^d)$, the $s$-proximal value is defined as:
\begin{equation}
\label{eqn: proximal-value}
    P_s(x) := \mathop{\arg\min}_{y\in\mr^d}\left\{ \frac{1}{2s}\left\| y - \left(x - s\nabla f(x)\right) \right\|^2 + g(y) \right\}, \quad \forall x \in \mathbb{R}^d.
\end{equation}
The $s$-proximal value $P_s(x)$ minimizes a weighted sum of the squared Euclidean distance from the gradient-based update $x - s\nabla f(x)$ and the regularization term $g(y)$. This formulation facilitates efficient minimization of composite objectives, even when $g(x)$ is non-smooth but proximally simple. For instance, when $g(x) = \lambda \|x\|_1$ (the $\ell_1$-norm), the $s$-proximal value admits a closed-form solution. Specifically, for any $x \in \mathbb{R}^d$, the $i$-th component of the $s$-proximal value $P_s(x)$ is given by: 
\[
P_s(x)_i = \big(\left|\left(x - s\nabla f(x)\right)_i\right| - \lambda s\big)_+ \text{sgn}\big(\left(x - s\nabla f(x)\right)_i\big), 
\]
where $i=1,\ldots,d$ and $\text{sgn}(\cdot)$ denotes the sign function. 
\end{defn}

%This closed-form expression is particularly useful, as it provides an efficient way to compute the $s$-proximal value for the $\ell_1$-regularization term, which is widely used in sparse optimization problems.  %Finally, we propose the definition of the $s$-proximal subgradient as follows.

\begin{defn}[$s$-Proximal Subgradient]
\label{defn: subgradient}
For any $x \in \mathbb{R}^d$, the $s$-proximal subgradient is defined as:
\begin{equation}
\label{eqn: subgradient}
G_s(x): = \frac{x - P_s(x)}{s},
\end{equation}
where the $s$-proximal value $P_s(x)$ is given by~\eqref{eqn: proximal-value}. The subgradient $G_s(x)$ generalizes the classical gradient by incorporating the proximal update with respect to the regularizer $g$, capturing the direction of improvement for the composite function $\Phi$.
\end{defn}
% based on the $s$-proximal update. It can be viewed as a generalization of the traditional gradient, incorporating the proximal step with respect to the regularizer $g$.

Using the $s$-proximal subgradient from~\Cref{defn: subgradient}, we extend the fundemental inequality~\eqref{eqn: fund-inq-smooth} to the composite setting. The proximal version is presented below, rigorously established in~\citep{li2024linear2}.  
\begin{lemma}[Lemma 4 in~\citep{li2024linear2}]
\label{lem: fund-inq-composite}
Let $\Phi = f + g$ be a composite function with $f \in \mathcal{S}_{\mu, L}^1(\mathbb{R}^d)$ and $g \in \mathcal{F}^0(\mathbb{R}^d)$. Then, the following inequality holds for any step size $s \in (0,1/L]$:
\begin{equation}
\label{eqn: fund-inq-composite}
\Phi(y - sG_s(y)) - \Phi(x) \leq \left\langle G_s(y), y - x \right\rangle - \underbrace{\frac{\mu}{2}\|y - x\|^2}_{\mathbf{SC}} - \frac{s}{2} \|G_s(y)\|^2, \quad \forall x,y \in \mathbb{R}^d.
\end{equation}
\end{lemma}
%Finally, we present a key inequality for the $s$-proximal subgradient, which will be used in the subsequent sections. This result, previously established in~\citep{li2024linear}, is stated as follows: 
%\begin{lemma}[Lemma 4.3 in~\citep{li2024linear}]
%\label{lem: key-inq-composite}
%Let $f \in \mathcal{S}_{\mu,L}^{1}(\mathbb{R}^d)$ and $g \in \mathcal{F}^0(\mathbb{R}^d)$. Then, the s-proximal subgradient, as defined in~\eqref{defn: subgradient}, satisfies the following inequality:
%\begin{equation}
%\label{eqn: key-inq-composite}
%\|G_s(y)\|^2 \geq 2\mu \left( \Phi(y - sG_s(y)) - \Phi(x^{\star}) \right),
%\end{equation}
%for any $ y \in \mathbb{R}^d$. 
%\end{lemma}

\section{Lyapnov analysis for forward-backward accelerated algorithms}
\label{sec: nag}

In this section, we first outline the process of constructing a novel Lyapunov function, following the principled way presented in~\citep{chen2022revisiting}.  Leveraging this approach, we then derive an $O(1/k^{2\alpha})$ convergence rate for the \texttt{NAG}-$\alpha$ method applied to strongly convex functions. A formal theorem is also provided to rigorously characterize this result.

%======================================================================%
\subsection{Construction of a novel Lyapunov function}
\label{subsec: construction}

To analyze the~\texttt{NAG}-$\alpha$ method, as given by~\eqref{eqn: nag-alpha-gradient} and~\eqref{eqn: nag-alpha-momentum}, we first reformulate it in the phase-space representation. This reformulation leverages the implicit-velocity scheme, extending approaches previously applied to the classical~\texttt{NAG} method~\citep{chen2022revisiting, li2024linear2}. To facilitate this process, we  introduce a new iterative sequence $\sqrt{s} v_{k} = x_{k} - x_{k-1}$, which transforms the algorithm into the following phase-space representation:
\begin{subequations}
\label{eqn: nag-alpha-phase}
    \begin{empheq}[left=\empheqlbrace]{align}
         & x_{k+1} - x_{k} = \sqrt{s}v_{k+1},                                                                                                                                           \label{eqn: nag-alpha-phase-x}   \\
         & v_{k+1} - v_{k} = - \left(1 - \frac{(k-1)^{\alpha}}{k^{\alpha}+rk^{\alpha-1}} \right)v_{k} - \sqrt{s} \nabla f(y_{k}),                     \label{eqn: nag-alpha-phase-v} 
    \end{empheq}
\end{subequations}
where the iterative sequence $\{y_k\}_{k=0}^{\infty}$ satisfies the following iterative relation:
\begin{equation}
\label{eqn: nag-alpha-x-y-v}
y_{k} = x_{k} + \frac{(k-1)^{\alpha}}{k^{\alpha}+rk^{\alpha-1}} \cdot \sqrt{s}v_{k}. 
\end{equation}

\begin{itemize}
\item[(\textbf{I})] \textbf{Construction of mixed energy} In the context of Lyapunov analysis, the focus is often shifted from the iterative point $x_{k}$ to the difference between the current iterate and the optimal solution $x^{\star}$. The difference, $x_k - x^{\star}$, provides  deeper insights into the convergence behavior of the algorithm. Reformulating~\eqref{eqn: nag-alpha-phase-x} in terms of this difference yields:
\begin{equation}
\label{eqn: nag-alpha-phase-x-star}
\left( x_{k+1} - x^{\star} \right) - \left( x_{k} - x^{\star} \right) = \sqrt{s}v_{k+1}
\end{equation}
Substituting this expression~\eqref{eqn: nag-alpha-phase-x-star} into~\eqref{eqn: nag-alpha-phase-v}, we derive the following iterative relation: 
\begin{multline}
\left[ \sqrt{s} k^{\alpha}v_{k+1} + rk^{\alpha-1}(x_{k+1} - x^{\star}) \right] - \left[ \sqrt{s} (k-1)^{\alpha}v_{k} + r(k-1)^{\alpha-1}(x_{k} - x^{\star}) \right]                                                      \\
= r \left[ k^{\alpha-1} - (k-1)^{\alpha-1} \right] (x_k - x^{\star}) - s \left( k^{\alpha} + rk^{\alpha-1} \right) \nabla f(y_{k}).                                       \label{eqn: iterative-nag-alpha-phase}
\end{multline}
This iterative relation leads naturally to the mixed energy, which plays a central role in our Lyapunov analysis. The mixed energy is defined as:
\begin{equation}
\label{eqn: mix-nag-alpha}
\mathcal{E}_{mix}(k) = \frac12 \left\| \sqrt{s} (k-1)^{\alpha}v_{k} + r(k-1)^{\alpha-1}(x_{k} - x^{\star})  \right\|^2.
\end{equation}
Next, we examine the change in mixed energy from one iteration to the next. Using the iterative relation~\eqref{eqn: iterative-nag-alpha-phase},  the iterative difference in mixed energy is given by:
\begin{align}
\mathcal{E}_{mix}&(k +1) - \mathcal{E}_{mix}(k)                                                                                                                                                                \nonumber \\
                         = & \bigg \langle r \left[ k^{\alpha-1} - (k-1)^{\alpha-1} \right](x_k - x^{\star})  - s \left( k^{\alpha} + rk^{\alpha-1} \right) \nabla f(y_{k}),    \nonumber \\
                            & \quad \sqrt{s} (k-1)^{\alpha}v_{k} +  \frac{ r }{2} \left[ k^{\alpha-1} + (k-1)^{\alpha-1} \right] (x_{k} - x^{\star}) - \frac{s }{2} \left( k^{\alpha} + rk^{\alpha-1} \right) \nabla f(y_{k})\bigg \rangle.                                                                                                                                                                                                                       \label{eqn: iter-diff-mix-nag-alpha-1}
\end{align}
Substituting the iterative relation~\eqref{eqn: nag-alpha-x-y-v} into the iterative difference~\eqref{eqn: iter-diff-mix-nag-alpha-1}, we obtain:
\begin{align}
\mathcal{E}_{mix}(k +1) \;- &\; \mathcal{E}_{mix}(k) \nonumber \\
                                     =  &\;\bigg \langle -  \frac{r \left[ k^{\alpha-1} - (k-1)^{\alpha-1} \right] }{k^{\alpha} + rk^{\alpha-1}} \cdot \sqrt{s} (k-1)^{\alpha} v_k \nonumber \\
                                         & \quad + r \left[ k^{\alpha-1} - (k-1)^{\alpha-1} \right](y_k - x^{\star})- s \left( k^{\alpha} + rk^{\alpha-1} \right) \nabla f(y_{k}),  \nonumber\\
                                         & \quad \left( 1 - \frac{r}{2} \cdot \frac{  k^{\alpha-1} + (k-1)^{\alpha-1} }{ k^{\alpha} + rk^{\alpha-1}  } \right) \sqrt{s} (k-1)^{\alpha}v_{k} \nonumber \\
                                         & \quad + \frac{r}{2} \left[ k^{\alpha-1} + (k-1)^{\alpha-1} \right]  (y_{k} - x^{\star}) - \frac{s \left( k^{\alpha} + rk^{\alpha-1} \right)}{2}  \nabla f(y_{k})\bigg \rangle. \label{eqn: iter-diff-mix-nag-alpha-2}
\end{align}
This difference involves several terms, which can be simplified further by introducing new notations:
\begin{subequations}
\begin{empheq}[left=\empheqlbrace]{align}
& A_1(k, \alpha) = k^{\alpha-1} - (k-1)^{\alpha-1},     \label{eqn: a1}   \\
& A_2(k, \alpha) = k^{\alpha-1} + (k-1)^{\alpha-1},    \label{eqn: a2}  \\
& B(k, \alpha, r) = k^{\alpha} + rk^{\alpha-1}.            \label{eqn: b}
\end{empheq}
\end{subequations}
Using these new definitions,~\eqref{eqn: a1}---\eqref{eqn: b}, the iterative difference~\eqref{eqn: iter-diff-mix-nag-alpha-2} in mixed energy becomes:
\begin{align}
 \mathcal{E}_{mix}(k +1) - \mathcal{E}_{mix}(k)     = & - \underbrace{s\sqrt{s}k^{\alpha} (k-1)^{\alpha} \left\langle \nabla f(y_k), v_k \right\rangle}_{\mathbf{I}}                \nonumber     \\
                                                                                   & -  rsk^{\alpha-1} B(k, \alpha, r)  \left\langle \nabla f(y_k), y_k - x^{\star} \right\rangle                                              \nonumber    \\
                                                                                   & + \frac{s^2 B^2(k, \alpha, r)}{2} \left\| \nabla f(y_k) \right\|^2                                                                                    \nonumber    \\
                                                                                   &- \frac{rs (k-1)^{2\alpha} A_1(k,\alpha) \left[  2k^{\alpha} + rA_1(k, \alpha) \right] }{ 2B^2(k, \alpha, r) } \|v_k\|^2      \nonumber     \\   
                                                                                   & + \frac{r^2A_1(k,\alpha)A_2(k,\alpha)}{2} \left\| y_k - x^{\star} \right\|^2                                                                  \nonumber     \\ 
                                                                                   & +   r\sqrt{s}(k-1)^{\alpha} A_1(k,\alpha) \left( 1- \frac{r A_{2}(k, \alpha)}{B(k, \alpha, r)} \right) \left\langle y_k - x^{\star}, v_k \right\rangle.                                                \label{eqn: iter-diff-mix-nag-alpha-3}
\end{align}
This simplified representation enables us to bound and analyze the convergence properties of the algorithm effectively.

\item[(\textbf{II})] \textbf{Construction of potential energy} We now proceed to define the potential energy, incorporating a dynamic coefficient $\tau(k)$, which is initially left undetermined: 
\begin{equation}\label{eqn: pot-nag-alpha}
    \mathcal{E}_{pot}(k) = s \tau(k) (f(x_k) - f(x^{\star})).
\end{equation}
Next, we compute the iterative difference in potential energy between consecutive iterations:
\begin{align}
 \mathcal{E}_{pot}(k+1)  - \mathcal{E}_{pot}(k)  = & s \tau(k) \left( f(x_{k+1}) - f(x_{k}) \right)                                                   \nonumber  \\
                                                                              & + s \left( \tau(k+1) - \tau(k) \right) \left( f(x_{k+1}) - f(x^{\star}) \right).       \label{eqn: iter-diff-pot-nag-alpha-1}
\end{align}
To bound this difference, we apply the fundamental inequality~\eqref{eqn: fund-inq-smooth} at two different points.  First, for $x_{k+1}$ and $x_{k}$,  we obtain:
\begin{equation}
 \label{eqn: fund-inq-xk-xk}
 f(x_{k+1}) - f(x_{k}) \leq \langle \nabla f(y_{k}), y_{k} - x_{k} \rangle - \underbrace{\frac{\mu}{2} \|y_{k} - x_{k}\|^2}_{\mathbf{SC}_1} - \frac{s}{2} \|\nabla f(y_{k})\|^2.  
\end{equation}
Similarly, applying the fundamental inequality~\eqref{eqn: fund-inq-smooth} to $x_{k+1}$ and $x^{\star}$, we derive:
\begin{equation}
 \label{eqn: fund-inq-xk-xstar}
 f(x_{k+1}) - f(x^{\star}) \leq \langle \nabla f(y_{k}), y_{k} - x^{\star} \rangle - \underbrace{\frac{\mu}{2} \|y_{k} - x^{\star}\|^2}_{\mathbf{SC}_2} - \frac{s}{2} \|\nabla f(y_{k})\|^2.  
\end{equation}
Here, the terms denoted $\mathbf{SC}_1$ in~\eqref{eqn: fund-inq-xk-xk} and $\mathbf{SC}_2$ in~\eqref{eqn: fund-inq-xk-xstar} are specific to the properties of strongly convex functions. These terms vanish in the case of convex functions, which aids in deriving convergence rates.  Substituting the two inequalities,~\eqref{eqn: fund-inq-xk-xk} and~\eqref{eqn: fund-inq-xk-xstar}, into the iterative difference~\eqref{eqn: iter-diff-pot-nag-alpha-1}, we establish the following upper bound:
\begin{align}
 \mathcal{E}_{pot}(k+1) - \mathcal{E}_{pot}(k) & \leq s \tau(k) \bigg( \left\langle \nabla f(y_k), y_k - x_k \right\rangle - \underbrace{ \frac{\mu}{2} \left\| y_k - x_k \right\|^2  }_{\mathbf{SC}_1}  \bigg)             \nonumber          \\
& \mathrel{\phantom{\leq}} + s \left( \tau(k+1) - \tau(k) \right)  \bigg( \left\langle \nabla f(y_k), y_k - x^{\star} \right\rangle - \underbrace{ \frac{\mu}{2} \left\| y_k - x^{\star} \right\|^2}_{\mathbf{SC}_2 }  \bigg) \nonumber          \\
& \mathrel{\phantom{\leq}} - \frac{s^2\tau(k+1)}{2} \cdot  \|\nabla f(y_k)\|^2          \label{eqn: iter-diff-pot-nag-alpha-2}
\end{align}
To simplify further, we substitute the iterative relation~\eqref{eqn: nag-alpha-x-y-v} into the right-hand side of~\eqref{eqn: iter-diff-pot-nag-alpha-2}. This leads to the following expression for the iterative difference in potential energy:
\begin{align}
 \mathcal{E}_{pot}(k+1) - \mathcal{E}_{pot}(k) \leq 
& \underbrace{s\sqrt{s} \tau(k) \cdot \frac{(k-1)^{\alpha}}{k^{\alpha}+rk^{\alpha-1}}  \cdot  \left\langle \nabla f(y_k), v_k \right\rangle}_{\mathbf{II}}                                                     \nonumber \\
& - \underbrace{\frac{\mu s^2}{2} \cdot \frac{\tau(k)(k-1)^{2\alpha}}{B^2(k,\alpha,r)} \left\| v_k \right\|^2}_{\mathbf{SC}_1}                                                                                              \nonumber \\
& + s \left( \tau(k+1) - \tau(k) \right)  \bigg( \left\langle \nabla f(y_k), y_k - x^{\star} \right\rangle - \underbrace{\frac{\mu}{2} \left\| y_k - x^{\star} \right\|^2}_{\mathbf{SC}_2} \bigg)    \nonumber \\
& - \frac{s^2\tau(k+1)}{2} \cdot  \| \nabla f(y_k) \|^2.                                                                                                                                                                  \label{eqn: iter-diff-pot-nag-alpha-3}
\end{align}
To ensure that the term $\mathbf{I} - \mathbf{II} = 0$,  we choose the dynamic coefficient $\tau(k)$ as: 
\begin{equation}
\label{eqn: pot-coefficient}
\tau(k) = k^\alpha B(k, \alpha, r)= k^\alpha \left( k^\alpha + rk^{\alpha-1} \right).
\end{equation}
Substituting the dynamic coefficient~\eqref{eqn: pot-coefficient} into~\eqref{eqn: iter-diff-pot-nag-alpha-3}, we simplify the upper bound, yielding the following expression:
\begin{align}
 \mathcal{E}_{pot}(k+1) - \mathcal{E}_{pot}(k) \leq & \underbrace{s\sqrt{s} k^{\alpha}(k-1)^{\alpha} \cdot  \left\langle \nabla f(y_k), v_k \right\rangle}_{\mathbf{II}}                             \nonumber \\
                                                                                & + s \left[ (k+1)^{\alpha}B(k+1, \alpha, r) - k^{\alpha}B(k, \alpha, r)  \right] \left\langle \nabla f(y_k), y_k - x^{\star} \right\rangle   \nonumber \\
                                                                                & - \frac{s^2(k+1)^{\alpha}B(k+1, \alpha, r)}{2}  \| \nabla f(y_k) \|^2                                                                                        \nonumber \\
                                                                                & - \underbrace{\frac{\mu s^2}{2} \cdot \frac{k^{\alpha}(k-1)^{2\alpha}}{B(k,\alpha, r)} \left\| v_k \right\|^2}_{\mathbf{SC}_1}            \nonumber \\
                                                                                & - \underbrace{\frac{\mu s}{2} \left[ (k+1)^{\alpha}B(k+1, \alpha, r) - k^{\alpha}B(k, \alpha, r)  \right]    \left\| y_k - x^{\star} \right\|^2}_{\mathbf{SC}_2}.               \label{eqn: iter-diff-pot-nag-alpha-4}
\end{align}
\end{itemize}

Finally, with the mixed energy defined in~\eqref{eqn: mix-nag-alpha} and the potential energy defined in~\eqref{eqn: pot-nag-alpha}, we can now construct the Lyapunov function as:
\begin{equation}
\label{eqn: lyapunov-nag-alpha}
\mathcal{E}(k) = sk^\alpha \left( k^\alpha + rk^{\alpha-1} \right) \left( f(x_k) - f(x^{\star}) \right) +  \frac12 \left\| \sqrt{s} (k-1)^{\alpha}v_{k} + r(k-1)^{\alpha-1}(x_{k} - x^{\star})  \right\|^2.
\end{equation}
This Lyapunov function encapsulates both the mixed and potential energies, serving as a key tool in deriving convergence rates for the~\texttt{NAG}-$\alpha$ method, as described in~\eqref{eqn: nag-alpha-gradient} and~\eqref{eqn: nag-alpha-momentum}.

%                              = & - s\sqrt{s}(k-1)^{\alpha}(k^{\alpha} + rk^{\alpha-1}) \left\langle \nabla f(y_k), v_k \right\rangle                                 \nonumber  \\
%                                 & + r\sqrt{s}(k-1)^{\alpha} \left[ k^{\alpha-1} - (k-1)^{\alpha-1} \right] \left\langle x_k - x^{\star}, v_k \right\rangle        \nonumber  \\
%                                 & + \frac{r^2 \left[ k^{2(\alpha-1)} - (k-1)^{2(\alpha-1)} \right]}{2} \cdot \|x_k - x^{\star}\|^2                                           \nonumber  \\
%                                 & - rsk^{\alpha-1} \left( k^{\alpha} + rk^{\alpha-1}\right) \left\langle \nabla f(y_k), x_k - x^{\star} \right\rangle              \nonumber  \\
%                                 & + \frac{s^2(k^{\alpha} + rk^{\alpha-1})^2}{2} \cdot \left\| \nabla f(y_k) \right\|^2                                                         \label{eqn: iter-diff-mix-nag-alpha-1} 

%======================================================================%
\subsection{Controllable $O\left( 1/k^{2\alpha} \right)$ convergence of NAG-$\alpha$}
\label{subsec: nag}
By leveraging the Lyapunov function defined in~\eqref{eqn: lyapunov-nag-alpha}, we rigorously establish the convergence rate for the~\texttt{NAG}-$\alpha$ method. This result is encapsulated  in the following theorem:

\begin{theorem}
\label{thm: nag-alpha}
Let $f \in \mathcal{S}_{\mu, L}^1(\mathbb{R}^d)$, and suppose the step size satisfies $s \in (0, 1/L]$. If the controllable parameter is chosen such that $r > 2\alpha$, then there exists a positive integer $K:= K(\alpha, r)$ such that the iterative sequence $\{x_k\}_{k=0}^{\infty}$, generated by \texttt{NAG}-$\alpha$ under updates~\eqref{eqn: nag-alpha-gradient} and~\eqref{eqn: nag-alpha-momentum} with any initial $x_0 = y_0 \in \mathbb{R}^d$, converges at the following rate:
\begin{equation}
\label{eqn: nag-alpha-rate}
f(x_k) - f(x^{\star}) \leq \frac{\mathcal{E}(K)}{sk^\alpha \left( k^\alpha + rk^{\alpha-1} \right)}, \quad \forall k \geq K.
\end{equation}
\end{theorem}

\begin{proof}[Proof of~\Cref{thm: nag-alpha}]
To establish the stated result, we analyze the iterative change in the Lyapunov function $\mathcal{E}(k)$ across iterations. By combining the inequalities derived for the iterative differences,~\eqref{eqn: iter-diff-mix-nag-alpha-3} and~\eqref{eqn: iter-diff-pot-nag-alpha-4}, we observe that the change in the Lyapunov function~\eqref{eqn: lyapunov-nag-alpha} satisfies:
\begin{align}
\mathcal{E}(k+1) - \mathcal{E}(k) \leq  &  \underbrace{s \left[ (k+1)^\alpha B(k+1, \alpha, r) - B^2(k, \alpha, r) \right]  \left( \left\langle \nabla f(y_k), y_k - x^{\star} \right\rangle - \frac{s}{2} \left\| \nabla f(y_k) \right\|^2  \right)}_{\mathbf{III}}                                                                                                                                                                                                                                        \nonumber     \\
                                                             & - \underbrace{\frac{rs (k-1)^{2\alpha} A_1(k,\alpha) \left[  2k^{\alpha} + rA_1(k, \alpha) \right] }{ 2B^2(k, \alpha, r) } \|v_k\|^2}_{\mathbf{IV}_1}     \nonumber     \\    
                                                              & + \underbrace{\frac{r^2A_1(k,\alpha)A_2(k,\alpha)}{2} \left\| y_k - x^{\star} \right\|^2}_{\mathbf{IV}_2}                                                                   \nonumber     \\ 
                                                              & + \underbrace{ r\sqrt{s}(k-1)^{\alpha} A_1(k,\alpha)  \left( 1- \frac{r A_{2}(k, \alpha)}{B(k, \alpha, r)} \right)  \left\langle y_k - x^{\star}, v_k \right\rangle}_{\mathbf{IV}_3}                                                   \nonumber     \\
                                                              & -  \underbrace{\frac{\mu s^2}{2} \cdot \frac{k^{\alpha}(k-1)^{2\alpha}}{B(k,\alpha, r)} \left\| v_k \right\|^2}_{\mathbf{SC}_1}                                          \nonumber     \\
                                                              & - \underbrace{\frac{\mu s}{2} \cdot \left[ (k+1)^{\alpha}B(k+1, \alpha, r) - k^{\alpha}B(k, \alpha, r) \right]  \left\| y_k - x^{\star} \right\|^2}_{\mathbf{SC}_2}. \label{eqn: lyapunov-iter-diff}   
\end{align}

To guarantee $\mathcal{E}(k+1) - \mathcal{E}(k) \leq 0$, we systematically analyze the contributions of each term by grouping them into four major components:

\begin{itemize}
\item[(1)]  \textbf{Analyzing Term $\mathbf{III}$.}  The coefficient of Term $\mathbf{III}$, after an asymptotic expansion, is given by:
\begin{align}
(k+1)^\alpha B(k+1, & \alpha, r)  - B^2(k, \alpha, r)                                                                                                                                                                                     \nonumber \\
                    & = (k+1)^{2\alpha} + r(k+1)^{2\alpha-1} - \left( k^{\alpha} + rk^{\alpha-1} \right)^2                                                                                                   \nonumber \\
                    & = k^{2\alpha} + \left( r + 2\alpha \right) k^{2\alpha-1} + O( k^{2 (\alpha - 1)} ) - k^{2\alpha} - 2r k^{2\alpha - 1} - r^2 k^{2\left( \alpha - 1 \right)}   \nonumber \\
                    & = \left( 2\alpha - r \right) k^{2\alpha-1} + O(k^{2 (\alpha - 1)} ).                                                                                                                              \label{eqn: iii-coefficient}  
\end{align}
It follows that when we control the parameter satisfing $r > 2 \alpha$, the coefficient becomes asymptotically non-positive, ensuring a negative contribution to the iterative difference. On the other hand, it follows from the fact that
\begin{equation*}
    \left\langle \nabla f(x) - \nabla f(y), x - y \right\rangle \geq \frac{1}{L} \left\| \nabla f(x) - \nabla f(y) \right\|^2, \quad \forall x,y \in \mathbb{R}^d
\end{equation*}
for $ f \in \mathcal{F}^1_L(\mathbb{R}^d) $ that 
\begin{equation*}
    \left\langle \nabla f(y_k), y_k - x^{\star} \right\rangle - \frac{s}{2} \left\| \nabla f(y_k) \right\|^2 \geq \left( \frac{1}{L} - \frac{s}{2} \right) \left\| \nabla f(y_k) \right\|^2 \geq 0.
\end{equation*}
Thus, the whole term $\mathbf{III}$ is negative in the iterative difference.

%where we observe that $(k+1)B(k+1, \alpha, r)  - B^2(k, \alpha, r)  \leq 0 $ when the controllable parameter $r > 2 \alpha$.

\item[(2)]  \textbf{Analyzing Term $\mathbf{IV}_1$.} The asymptotic expansion of the coefficient of Term $\mathbf{IV}_1$ is:%\begin{equation}%\label{eqn: iv-11}%\end{equation}
\[
- \frac{rs (k-1)^{2\alpha} A_1(k,\alpha) \left[  2k^{\alpha} + rA_1(k, \alpha) \right] }{ 2B^2(k, \alpha, r) } = - rs(\alpha-1) k^{2\alpha-2} + O(k^{2\alpha-3}).
\]
For $\alpha < 1$, this term is asymptotically non-negative. Incorporating the term $\mathbf{SC}_1$ from the strongly convex inequality, we obtain: 
\begin{equation}
\label{eqn: iv-12}
 - \frac{rs (k-1)^{2\alpha} A_1(k,\alpha) \left[  2k^{\alpha} + rA_1(k, \alpha) \right] }{ 2B^2(k, \alpha, r) } - \frac{\mu s^2}{2} \cdot \frac{k^{\alpha}(k-1)^{2\alpha}}{B(k,\alpha, r)}  =  - \frac{\mu s^2}{2} \cdot k^{2\alpha} + O(k^{2\alpha-1}).
\end{equation}
Thus, the contribution of $\|v_k\|^2$ remains asymptotically negative. 

\item[(3)]  \textbf{Analyzing Term $\mathbf{IV}_2$.}  The asymptotic expansion of the coefficient of Term $\mathbf{IV}_2$ is:%\end{equation}%\begin{equation}%\label{eqn: iv-21}
\[
\frac{r^2A_1(k,\alpha)A_2(k,\alpha)}{2} = r^2(\alpha - 1)k^{2\alpha - 3} + O(k^{2\alpha-4}).
\]
For $\alpha > 1$, this term is asymptotically non-negative. Incorporating the term $\mathbf{SC}_2$ from the strongly convex inequality, we obtain: 
\begin{multline}
\label{eqn: iv-22}
\frac{r^2A_1(k,\alpha)A_2(k,\alpha)}{2}  - \frac{\mu s}{2} \cdot \left[ (k+1)^{\alpha}B(k+1, \alpha, r) - k^{\alpha}B(k, \alpha, r) \right] \\ = - \mu s \alpha k^{2\alpha-1} + O(k^{2\alpha-2}).
\end{multline}
Thus, the contribution of $\|y_k - x^{\star}\|^2$ also remains asymptotically negative. 

\item[(4)]  \textbf{Analyzing Term $\mathbf{IV}_3$.}  The coefficient of Term $\mathbf{IV}_3$ is asymptotically expanded as:
\begin{equation}
\label{eqn: iv-3}
r\sqrt{s}(k-1)^{\alpha} A_1(k,\alpha) \left( 1- \frac{r A_{2}(k, \alpha)}{B(k, \alpha, r)} \right)  = r\sqrt{s}(\alpha - 1)k^{2\alpha-2} + O( k^{2\alpha - 3}).
\end{equation}
By combining the three asymptotic expansions,~\eqref{eqn: iv-12},~\eqref{eqn: iv-22}, and~\eqref{eqn: iv-3}, the two strongly convex terms, $\mathbf{SC}_1$ and $\mathbf{SC}_2$, reduce the combined sum $-\mathbf{IV}_1 + \mathbf{IV}_2 + \mathbf{IV}_3$ to a sum of three non-positive perfect squares. %Together with the inequality~\eqref{eqn: iii-coefficient}, the overall iterative difference is guaranteed to be non-positive, thereby ensuring $\mathcal{E}(k+1) - \mathcal{E}(k) \leq 0$.
\end{itemize}
In summary, we conclude that when $r > 2\alpha$, there exists a positive integer $K=K(\alpha, r)$ such that for any $k \geq K$, the iterative difference $\mathcal{E}(k+1) - \mathcal{E}(k)$ remains non-positive.  This ensures that the Lyapunov function~\eqref{eqn: lyapunov-nag-alpha} is non-increasing, thereby completing the proof.

\end{proof}

From the above proof, it becomes clear that the fundamental inequality~\eqref{eqn: fund-inq-smooth} for strongly convex functions is pivotal in the analysis and cannot be simplified to its counterpart for convex functions. This distinction is further elaborated in the following remark:
\begin{remark}
In the inequality governing the iterative difference~\eqref{eqn: lyapunov-iter-diff}, the coefficients of both terms, $\mathbf{IV}_1$ and $\mathbf{IV}_2$, involve the expression 
\[
A_1(k,\alpha) = k^{\alpha - 1} - (k-1)^{\alpha - 1} = (\alpha - 1)k^{\alpha - 2} + O\left( k^{\alpha  - 3}\right),
\]
which retains a consistent sign based on the value of $\alpha$. Specifically, when $\alpha \geq 1$, the sign is positive, while for $\alpha < 1$, the sign is negative. As a result,  in the inequality for the iterative difference~\eqref{eqn: lyapunov-iter-diff}, the two terms, $-\mathbf{IV}_1$ and $\mathbf{IV}_2$, always exhibits opposite signs. This observation reveals an important limitation: if the objective function under consideration is merely convex, rather than strongly convex, the terms $\mathbf{SC}_1$ and $\mathbf{SC}_2$, which arise from the fundamental inequaliy for strongly convex functions,  vanish. Without these additional terms, it is no longer possible to ensure that the iterative difference~\eqref{eqn: lyapunov-iter-diff} is non-positive. Consequently, the Lyapunov function~\eqref{eqn: lyapunov-nag-alpha} cannot be guaranteed to decrease consistently. This underscores the indispensable role of strong convexity in achieving the desired convergence properties. The additional structure provided by strong convexity ensures that the iterative process maintains the non-positivity of the iterative difference, which is essential for guaranteeing a monotonic decrease in the Lyapunov function.
\end{remark}

With~\Cref{thm: nag-alpha}, we can also derive the convergence rate of~\texttt{NAG}-$\alpha$ at the critical step size $s = 1/L$, filling a gap left in prior analyses such as those in~\citep{li2024linear, fu2024lyapunov}. 
\begin{corollary}
\label{coro: nag-alpha-l}
Under the same assumptions of~\Cref{thm: nag-alpha}, if the step size satisfies $s = 1/L$, there exists a positive integer $K:= K(\alpha, r)$ such that \texttt{NAG}-$\alpha$ converges at the following rate:
\begin{equation}
\label{eqn: nag-alpha-rate-l}
f(x_k) - f(x^{\star}) \leq \frac{L\mathcal{E}(K)}{k^\alpha \left( k^\alpha + rk^{\alpha-1} \right)}, \quad \forall k \geq K.
\end{equation}
\end{corollary}
This result confirms the theoretical guarantees of \texttt{NAG}-$\alpha$ for strongly convex functions and provides explicit conditions for achieving optimal convergence rates under the critical step size.

%======================================================================%
\subsection{Generalization to FISTA-$\alpha$}
\label{subsec: fista}

We now extend the controllable $O\left( 1/k^{2\alpha} \right)$ convergence of~\texttt{NAG}-$\alpha$, as established in~\Cref{thm: nag-alpha} and~\Cref{coro: nag-alpha-l}, to its proximal counterpart,~\texttt{FISTA}-$\alpha$. In accordance with~\Cref{defn: proximal-value},~\texttt{FISTA}-$\alpha$ employs the $s$-proximal value~\eqref{eqn: proximal-value} and adheres to the following iterative scheme, starting from any initial point $y_0 = x_0 \in \mathbb{R}^d$: 
\begin{subequations}
\label{eqn: fista-alpha}
    \begin{empheq}[left=\empheqlbrace]{align}
         & x_{k} = P_s(y_{k-1})                                                                                                                  \label{eqn: fista-alpha-gradient}         \\
         & y_{k} = x_{k} + \frac{(k-1)^{\alpha}}{k^{\alpha}+rk^{\alpha-1}} \left( x_{k} - x_{k-1} \right)         \label{eqn: fista-alpha-momentum} 
    \end{empheq}
\end{subequations}
where $s>0$ is the step size. As specified in~\Cref{defn: subgradient},~\texttt{FISTA}-$\alpha$ can alternatively be expressed using the $s$-proximal subgradient~\eqref{eqn: subgradient}, yielding the alternative iterative scheme:
\begin{subequations}
\label{eqn: fista-alpha}
    \begin{empheq}[left=\empheqlbrace]{align}
         & x_{k} = y_{k-1} - sG_s(y_{k-1})                                                                                                  \label{eqn: fista1-alpha-gradient}         \\
         & y_{k} = x_{k} + \frac{(k-1)^{\alpha}}{k^{\alpha}+rk^{\alpha-1}} \left( x_{k} - x_{k-1} \right)         \label{eqn: fista1-alpha-momentum} 
    \end{empheq}
\end{subequations}
where the $s$-proximal subgradient $G_s(y_k)$ generalizes the classical gradient $\nabla f(y_k)$ used in~\texttt{NAG}-$\alpha$. This formulation underscores the structural similarity between \texttt{NAG}-$\alpha$ and \texttt{FISTA}-$\alpha$, while explicitly incorporating the non-smooth nature of the composite objective function $\Phi$. By introducing the velocity iterative sequence $v_k = \left( x_k - x_{k-1} \right)/\sqrt{s}$, the~\texttt{FISTA}-$\alpha$ updates can be formulated in a phase-space representation akin to \texttt{NAG}-$\alpha$:
\begin{subequations}
\label{eqn: fista-alpha-phase}
    \begin{empheq}[left=\empheqlbrace]{align}
         & x_{k+1} - x_{k} = \sqrt{s}v_{k+1},                                                                                                                                           \label{eqn: fista-alpha-phase-x}   \\
         & v_{k+1} - v_{k} = - \left(1 - \frac{(k-1)^{\alpha}}{k^{\alpha}+rk^{\alpha-1}} \right)v_{k} - \sqrt{s} G_s(y_{k}),                          \label{eqn: fista-alpha-phase-v} 
    \end{empheq}
\end{subequations}
with the auxiliary sequence $\{y_k\}_{k=0}^{\infty}$ satisfying:
\begin{equation}
\label{eqn: fista-alpha-x-y-v}
y_{k} = x_{k} + \frac{(k-1)^{\alpha}}{k^{\alpha}+rk^{\alpha-1}} \cdot \sqrt{s}v_{k}. 
\end{equation}

To establish the controllable $O\left( 1/k^{2\alpha} \right)$ convergence rate for~\texttt{FISTA}-$\alpha$, we construct a generalized Lyapunov function inspired by the smooth one~\eqref{eqn: lyapunov-nag-alpha}.   This Lyapunov function is adapted to composite functions by substituting $f$ with $\Phi = f + g$, yielding:
\begin{equation}
\label{eqn: lyapunov-fista-alpha}
\mathcal{E}(k) = sk^\alpha \left( k^\alpha + rk^{\alpha-1} \right) \left( \Phi(x_k) - \Phi(x^{\star}) \right) +  \frac12 \left\| \sqrt{s} (k-1)^{\alpha}v_{k} + r(k-1)^{\alpha-1}(x_{k} - x^{\star})  \right\|^2.
\end{equation}
In~\Cref{subsec: construction}, the $O\left( 1/k^{2\alpha} \right)$ rate for smooth functions hinges on two key inequalities,~\eqref{eqn: fund-inq-xk-xk} and~\eqref{eqn: fund-inq-xk-xstar}. For the composite setting, these are generalized to accommodate proximal components by leveraging~\Cref{lem: fund-inq-composite}, which extends the fundamental inequality for strongly convex functions to the proximal setting. These adaptations enable a rigorous analysis of~\texttt{FISTA}-$\alpha$, effectively addressing the challenges posed by the non-smooth components in $\Phi$. The convergence properties are formally encapsulated in the following theorem:
\begin{theorem}
\label{thm: fista-alpha}
Let $f \in \mathcal{S}_{\mu, L}^1(\mathbb{R}^d)$ and $g \in \mathcal{F}^{0}(\mathbb{R}^d)$, and suppose the step size satisfies $s \in (0, 1/L]$. If the controllable parameter is chosen such that $r > 2\alpha$, then there exists a positive integer $K:= K(\alpha, r)$ such that the iterative sequence $\{x_k\}_{k=0}^{\infty}$, generated by \texttt{FISTA}-$\alpha$ under updates~\eqref{eqn: fista-alpha-gradient} and~\eqref{eqn: fista-alpha-momentum} with any initial $x_0 = y_0 \in \mathbb{R}^d$, converges at the following rate:
\begin{equation}
\label{eqn: fista-alpha-rate}
\Phi(x_k) - \Phi(x^{\star}) \leq \frac{\mathcal{E}(K)}{sk^\alpha \left( k^\alpha + rk^{\alpha-1} \right)}, \quad \forall k \geq K.
\end{equation}
Furthermore, if the step size is set critically as $s = 1/L$, the iterative sequence $\{x_k\}_{k=0}^{\infty}$ satisfies:
\begin{equation}
\label{eqn: fista-alpha-rate-l}
\Phi(x_k) - \Phi(x^{\star}) \leq \frac{L\mathcal{E}(K)}{k^\alpha \left( k^\alpha + rk^{\alpha-1} \right)}, \quad \forall k \geq K.
\end{equation}
\end{theorem}

\section{Monotonically controllable $O\left( 1/k^{2\alpha} \right)$ convergence}
\label{sec: m-nag-fista}

In this section, we extend the controllable $O\left( 1/k^{2\alpha} \right)$ convergence rate established for \texttt{NAG}-$\alpha$ in~\Cref{sec: nag} to its proximal variant, \texttt{M-NAG}-$\alpha$. Following the methodology outlined by~\citet{fu2024lyapunov}, the Lyapunov function~\eqref{eqn: lyapunov-nag-alpha} is designed without incorporating the kinetic energy term. As a result, the derivation of  the convergence rate relies solely on the iterative dynamics of \texttt{M-NAG}-$\alpha$, as specified in~\eqref{eqn: m-nag-gradient} ---~\eqref{eqn: m-nag-momentum}, bypassing the need for the complete set of~\texttt{NAG} updates in~\eqref{eqn: nag-alpha-gradient} and~\eqref{eqn: nag-alpha-momentum}. Building on this framework, we further generalize the controllable $O\left( 1/k^{2\alpha} \right)$ convergence to its proximal counterpart,~\texttt{M-FISTA}-$\alpha$. 

%%%%%%%%%%%%%%%%%%%%%%%%%%%%%%%%%%%%%%%%%%%%%%%%%%%%%%%%%%%%%%%%%%%%%%%%%%%%%%%%%%%%%%%%%%%%
\subsection{Smooth optimization via M-NAG-$\alpha$}
\label{subsec: m-nag}

According to the iterative relation~\eqref{eqn: nag-alpha-x-y-v}, the velocity term $v_k$ can be expressed in terms of $x_k$ and $y_{k}$. This substitution yields the following relation:
\begin{multline}
\label{eqn: nag-to-m-nag}
 \sqrt{s} (k-1)^{\alpha}v_{k} + r(k-1)^{\alpha-1}(x_{k} - x^{\star}) \\ = \left[ k^{\alpha} + r\left( k^{\alpha-1} - (k-1)^{\alpha-1} \right) \right] \left( y_{k} - x_{k}\right)  +  r(k-1)^{\alpha - 1}\left( y_k - x^{\star} \right). 
\end{multline}
Leveraging this expression, we reformulate the Lyapunov function originally introduced in~\eqref{eqn: lyapunov-nag-alpha} so that it depends exclusively on the variables, $x_k$ and $y_k$. The reformulated Lyapunov function is given as:
\begin{align}
\mathcal{E}(k) = & \underbrace{sk^{\alpha} \left( k^{\alpha} + rk^{\alpha-1} \right)\left( f(x_{k}) - f(x^{\star}) \right)}_{:=\;\mathcal{E}_{pot}(k)}  \nonumber \\
                           &+ \underbrace{\frac12\left\| \left[ k^{\alpha} + r\left( k^{\alpha-1} - (k-1)^{\alpha-1} \right) \right] \left( y_{k} - x_{k}\right)  +  r(k-1)^{\alpha - 1}\left( y_k - x^{\star} \right) \right\|^2}_{:=\;\mathcal{E}_{mix}(k)}  \label{eqn: mnag-lyapunov}
\end{align}
This reformulation, which decouples $\mathcal{E}(k)$ from $v_k$ and depends only on $x_k$ and $y_k$, is instrumental in deriving the controllable  $O\left( 1/k^{2\alpha} \right)$ convergence rate for \texttt{M-NAG}-$\alpha$. The formal result is encapsulated in the subsequent theorem.
\begin{theorem}
\label{thm: m-nag-alpha}
Let $f \in \mathcal{S}_{\mu, L}^1(\mathbb{R}^d)$, and suppose the step size satisfies $s \in (0, 1/L]$. If the controllable parameter is chosen such that $r > 2\alpha$, then there exists a positive integer $K:= K(\alpha, r)$ such that the iterative sequence $\{x_k\}_{k=0}^{\infty}$, generated by \texttt{M-NAG}-$\alpha$ under updates~\eqref{eqn: m-nag-gradient} ---~\eqref{eqn: m-nag-momentum} with any initial $x_0 = y_0 \in \mathbb{R}^d$, converges monotonically at the following rate:
\begin{equation}
\label{eqn: m-nag-alpha-rate}
f(x_k) - f(x^{\star}) \leq \frac{\mathcal{E}(K)}{sk^\alpha \left( k^\alpha + rk^{\alpha-1} \right)}, \quad \forall k \geq K.
\end{equation}
Furthermore, if the step size is set critically as $s = 1/L$, the iterative sequence $\{x_k\}_{k=0}^{\infty}$ satisfies:
\begin{equation}
\label{eqn: m-nag-alpha-rate-l}
f(x_k) - f(x^{\star}) \leq \frac{L\mathcal{E}(K)}{k^\alpha \left( k^\alpha + rk^{\alpha-1} \right)}, \quad \forall k \geq K.
\end{equation}
\end{theorem}

\begin{proof}[Proof of~\Cref{thm: m-nag-alpha}]
To calculate the iterative difference of the Lyapunov function, we break the process into two distinct components: the potential energy and the mixed energy, as shown in~\eqref{eqn: mnag-lyapunov}. To further streamline the calculation, we utilize the new definitions introduced in~\eqref{eqn: a1}---\eqref{eqn: b}. These definitions provide convenient expressions  that facilitate the management of interactions between the various terms in the Lyapunov function.
\begin{itemize}
\item[(\textbf{1})] For the potential energy $\mathcal{E}_{pot}(k)$, we calculate its iterative difference as
\begin{align}
\mathcal{E}_{pot}(k+1) - \mathcal{E}_{pot}(k) = & s k^{\alpha}B(k, \alpha, r) \left( f(x_{k+1}) - f(x_{k}) \right)                                                                                   \nonumber        \\
                                                                           & + s \left[ (k+1)^{\alpha}B(k+1, \alpha, r) - k^{\alpha}B(k, \alpha, r) \right] \left( f(x_{k+1}) - f(x^{\star}) \right) .    \label{eqn: pot-iter-diff-mnag-1} 
\end{align}
By combining~\eqref{eqn: m-nag-gradient} and~\eqref{eqn: m-nag-comparison}, we can deduce that the iterative sequence of~\texttt{M-NAG}-$\alpha$ satisfies: 
\begin{equation}
\label{eqn: m-nag-special}
f(x_{k+1}) \leq f(z_k) = f(y_{k} - s \nabla f(y_{k})).
\end{equation}
Thus, we can apply the fundamental inequality~\eqref{eqn: fund-inq-smooth}. First, substituting the points $x_{k+1}$ and $x_{k}$,  we obtain the following inequality:
\begin{align}
 f(x_{k+1}) - f(x_{k}) &  \leq f(y_{k} - s \nabla f(y_{k})) - f(x_k)     \nonumber       \\
                                &  \leq \langle \nabla f(y_{k}), y_{k} - x_{k} \rangle - \underbrace{\frac{\mu}{2} \|y_{k} - x_{k}\|^2}_{\mathbf{SC}_1} - \frac{s}{2} \|\nabla f(y_{k})\|^2.   \label{eqn: fund-inq-xk-xk-mnag}
\end{align}
Then, for the points $x_{k+1}$ and $x^{\star}$, we derive:
\begin{align}
 f(x_{k+1}) - f(x^{\star}) & \leq f(y_{k} - s \nabla f(y_{k})) - f(x^{\star})     \nonumber       \\ 
                                     & \leq \langle \nabla f(y_{k}), y_{k} - x^{\star} \rangle - \underbrace{\frac{\mu}{2} \|y_{k} - x^{\star}\|^2}_{\mathbf{SC}_2} - \frac{s}{2} \|\nabla f(y_{k})\|^2.   \label{eqn: fund-inq-xk-xstar-mnag}
\end{align}
Substituting the two inequalities,~\eqref{eqn: fund-inq-xk-xk-mnag} and~\eqref{eqn: fund-inq-xk-xstar-mnag}, into~\eqref{eqn: pot-iter-diff-mnag-1}, we establish the following upper bound for the iterative difference in the potential energy as:
\begin{align}
\mathcal{E}_{pot}(k+1) - \mathcal{E}_{pot}(k) & = \underbrace{s k^{\alpha}B(k, \alpha, r)  \langle \nabla f(y_{k}), y_{k} - x_{k} \rangle}_{\mathbf{I}}                                   \nonumber   \\
                                                                        & \mathrel{\phantom{=}} + s \left[ (k+1)^{\alpha}B(k+1, \alpha, r) - k^{\alpha}B(k, \alpha, r) \right] \langle \nabla f(y_{k}), y_{k} - x^{\star} \rangle                                                                                                                                                                                                                                                                \nonumber   \\
                                                                        & \mathrel{\phantom{=}}  - \frac{s^2(k+1)^{\alpha}B(k+1, \alpha, r)}{2}  \|\nabla f(y_{k})\|^2                                                     \nonumber   \\
                                                                        & \mathrel{\phantom{=}}  - \underbrace{\frac{\mu s}{2} \cdot k^{\alpha}B(k, \alpha, r) \|y_{k} - x_{k}\|^2}_{\mathbf{SC}_1}      \nonumber   \\
                                                                        &  \mathrel{\phantom{=}} - \underbrace{\frac{\mu s}{2} \cdot \left[ (k+1)^{\alpha}B(k+1, \alpha, r) - k^{\alpha}B(k, \alpha, r) \right] \|y_{k} - x^{\star}\|^2}_{\mathbf{SC}_2}.       \label{eqn: pot-iter-diff-mnag-2} 
\end{align}

\item[(\textbf{2})] Next, we turn our attention to the mixed energy.  By substituting~\eqref{eqn: m-nag-gradient} into~\eqref{eqn: m-nag-momentum} and re-expressing the differences between the iterative sequences and the optimal points,  rather than using the sequences themselves, we arrive at the following iterative relation:
\begin{multline} 
\big[ \left( (k+1)^{\alpha} + rA_1(k+1, \alpha) \right) \left( y_{k+1} - x_{k+1} \right) +rk^{\alpha-1} \left( y_{k+1} - x^{\star}\right)  \big]                                                                        \\
- \left[ \left( k^{\alpha} + rA_1(k, \alpha) \right) \left( y_{k} - x_{k} \right) +r(k-1)^{\alpha-1} \left( y_{k} - x^{\star}\right)  \right]                    \mathrel{\phantom{=======}}                 \\
= - rA_1(k, \alpha)\left( y_k - x_k \right) + rA_1(k, \alpha)\left( y_k - x^{\star} \right) - sB(k, \alpha, r) \nabla f(y_k ).                                  \label{eqn: iterative-nag-alpha-m}
\end{multline}
Using the iterative relation~\eqref{eqn: iterative-nag-alpha-m}, we can compute the iterative difference in the mixed energy as follows:
\begin{align}
\mathcal{E}_{mix}(k+&1)  - \mathcal{E}_{mix}(k) \nonumber  \\
= & \; \bigg \langle - rA_1(k, \alpha)\left( y_k - x_k \right) + rA_1(k, \alpha)\left( y_k - x^{\star} \right) - sB(k, \alpha, r) \nabla f(y_k ),  \nonumber \\
   & \; \; \left( k^{\alpha} + \frac{rA_1(k, \alpha)}{2} \right) \left( y_k - x_k \right) +   \frac{rA_2(k, \alpha)}{2}  \left( y_{k} - x^{\star}\right)   - \frac{s B(k, \alpha, r)}{2}  \nabla f(y_k) \bigg \rangle  \label{eqn: mix-iter-diff-mnag-1} 
\end{align}
Reformulating and simplifying the iterative difference~\eqref{eqn: mix-iter-diff-mnag-1}, we have:
\begin{align}
\mathcal{E}_{mix}(k+1)  - \mathcal{E}_{mix}(k) = &  - \underbrace{sk^{\alpha} B(k, \alpha, r) \left\langle \nabla f(y_k), y_k - x_k \right\rangle}_{\mathbf{II}}                   \nonumber \\
                                                                              &  - rsk^{\alpha-1} B(k, \alpha, r) \left\langle \nabla f(y_k), y_k - x^{\star} \right\rangle                                                     \nonumber \\
                                                                              & + \frac{s^2B^2(k, \alpha, r)}{2} \| \nabla f(y_k) \|^2                                                                                                     \nonumber \\
                                                                              & - \frac{r A_1(k, \alpha) \left[ 2k^{\alpha} + rA_1(k, \alpha) \right]}{2} \| y_k - x_k\|^2                                                   \nonumber \\
                                                                              & + \frac{r^2A_1(k, \alpha)A_2(k, \alpha)}{2} \|y_k - x^{\star} \|^2                                                                                  \nonumber \\
                                                                              & + r A_1(k, \alpha)\left( B(k, \alpha, r) - r A_2(k, \alpha) \right) \left\langle y_k - x_k, y_k - x^{\star} \right\rangle.       \label{eqn: mix-iter-diff-mnag-2} 
\end{align}
\end{itemize}

From the iterative differences of the potential and the mixed energies, as described in~\eqref{eqn: pot-iter-diff-mnag-2} and~\eqref{eqn: mix-iter-diff-mnag-2}, we observe that the contributions from terms $\mathbf{I}$ and $\mathbf{II}$ are equal. This observation allows us to derive the following bound on the iterative difference:
\begin{align}
\mathcal{E}(k+1) - \mathcal{E}(k) \leq  &  \underbrace{s \left[ (k+1)^\alpha B(k+1, \alpha, r) - B^2(k, \alpha, r) \right]  \left( \left\langle \nabla f(y_k), y_k - x^{\star} \right\rangle - \frac{s}{2} \left\| \nabla f(y_k) \right\|^2  \right)}_{\mathbf{III}}                                                                                                                                                                                                                                        \nonumber     \\
                                                             & -   \underbrace{  \frac{r A_1(k, \alpha) \left[ 2k^{\alpha} + rA_1(k, \alpha) \right]}{2} \| y_k - x_k\|^2   }_{\mathbf{IV}_1}     \nonumber     \\    
                                                             & +  \underbrace{  \frac{r^2A_1(k,\alpha)A_2(k,\alpha)}{2} \left\| y_k - x^{\star} \right\|^2}_{\mathbf{IV}_2}                                                                   \nonumber     \\ 
                                                              & + \underbrace{ r A_1(k, \alpha)\left( B(k, \alpha, r) - r A_2(k, \alpha) \right) \left\langle y_k - x_k, y_k - x^{\star} \right\rangle   }_{\mathbf{IV}_3}                                                   \nonumber     \\
                                                              & -  \underbrace{ \frac{\mu s}{2} \cdot k^{\alpha}B(k, \alpha, r) \|y_{k} - x_{k}\|^2}_{\mathbf{SC}_1}                                          \nonumber     \\
                                                              & - \underbrace{\frac{\mu s}{2} \cdot \left[ (k+1)^{\alpha}B(k+1, \alpha, r) - k^{\alpha}B(k, \alpha, r) \right]  \left\| y_k - x^{\star} \right\|^2}_{\mathbf{SC}_2}. \label{eqn: mnag-lyapunov-iter-diff}   
\end{align}

Analogously, we can estimate the asymptotic expansion for the iterative difference as follows:
\begin{itemize}
\item[(\textbf{1})] Coefficient of Term $\mathbf{III}$:
\begin{equation}
\label{eqn: m-nag-iii}
(k+1)^\alpha B(k+1, \alpha, r) - B^2(k, \alpha, r) = (2\alpha - r) k^{2\alpha-1} + O(k^{2(\alpha - 1)}).
\end{equation}
\item[(\textbf{2})] Coefficient of $\|y_k - x_k\|^2$ (Term $\mathbf{IV}_1$ and Term $\mathbf{SC}_1$):
\begin{equation}
\label{eqn: m-nag-iv-1}
- \frac{r A_1(k, \alpha) \left[ 2k^{\alpha} + rA_1(k, \alpha) \right]}{2} - \frac{\mu s}{2} \cdot k^{\alpha}B(k, \alpha, r) = - \frac{\mu s}{2} \cdot k^{2\alpha} + O(k^{2\alpha - 1}).
\end{equation}
\item[(\textbf{3})] Coefficient of $\|y_k - x^{\star}\|^2$ (Term $\mathbf{IV}_2$ and Term $\mathbf{SC}_2$):
\begin{multline}
\label{eqn: m-nag-iv-2}
\frac{r A_1(k, \alpha) A_2(k, \alpha)}{2} - \frac{\mu s}{2} \left[ (k+1)^{\alpha}B(k+1, \alpha, r) - k^{\alpha}B(k, \alpha, r) \right] \\ = - \mu s \alpha k^{2\alpha-1} + O(k^{2\alpha-2}).
\end{multline}
\item[(\textbf{4})] Coefficient of $\left\langle y_k - x_k, y_k - x^{\star} \right\rangle$ (Term $\mathbf{IV}_3$):
\begin{equation}
\label{eqn: m-nag-iv-3}
r A_1(k, \alpha)\left( B(k, \alpha, r) - r A_2(k, \alpha) \right) = (\alpha - 1) r k^{2(\alpha - 1)} + O(k^{2\alpha - 3}).
\end{equation}
\end{itemize}
By plugging~\eqref{eqn: m-nag-iii} ---~\eqref{eqn: m-nag-iv-3} into the iterative inequality~\eqref{eqn: mnag-lyapunov-iter-diff}, we can deduce that $\mathcal{E}(k+1) - \mathcal{E}(k) \leq 0$ asymptotically. This confirms that the Lyapunov function decreases over iterations, leading to the conclusion that the proof is complete.
\end{proof}

%%%%%%%%%%%%%%%%%%%%%%%%%%%%%%%%%%%%%%%%%%%%%%%%%%%%%%%%%%%%%%%%%%%%%%%%%%%%%%%%%%%%%%%%%%%%
\subsection{Composite optimization via M-FISTA-$\alpha$}
\label{subsec: m-fista}

In this section, we extend the controllable $O\left( 1/k^{2\alpha} \right)$ convergence rate of~\texttt{M-NAG}-$\alpha$, as established in~\Cref{thm: m-nag-alpha}, to include its proximal variant,~\texttt{M-FISTA}-$\alpha$. As outlined in~\Cref{defn: proximal-value}, \texttt{M-FISTA}-$\alpha$ utilizes the $s$-proximal value defined in~\eqref{eqn: proximal-value} and adheres to the following iterative scheme outlined below, starting from any initial point $y_0 = x_0 \in \mathbb{R}^d$:
\begin{subequations}
\begin{empheq}[left=\empheqlbrace]{align}
& z_{k-1} = P_s(y_{k-1}),                                                                                                                                       \label{eqn: m-fista-gradient}       \\
& x_{k} = \left\{ \begin{aligned} 
                          & z_{k-1}, \quad \text{if}\; f(z_{k-1}) \leq f(x_{k-1}),  \\
                          & x_{k-1}, \quad \text{if}\; f(z_{k-1}) > f(x_{k-1}),
                          \end{aligned} \right.                                                                                                                    \label{eqn: m-fista-comparison}   \\    
& y_{k} = x_{k} + \frac{(k-1)^{\alpha}}{k^{\alpha}+rk^{\alpha-1}} \left( x_{k} - x_{k-1} \right) + \frac{(k-1)^{\alpha} + r(k-1)^{\alpha-1}}{k^{\alpha}+rk^{\alpha-1}} \left(z_{k-1} - x_{k} \right),            \label{eqn: m-fista-momentum}
\end{empheq}
\end{subequations}
where $s>0$ is the step size. To offer a more unified and efficient perspective, we reformulate~\texttt{M-FISTA}-$\alpha$, which mirrors the structure of \texttt{M-NAG}-$\alpha$, with one key difference: the $s$-proximal value is replaced with the $s$-proximal subgradient as defined in~\eqref{defn: subgradient}. This adjustment results in the following iterative scheme:
\begin{subequations}
\begin{empheq}[left=\empheqlbrace]{align}
& z_{k-1} = y_{k-1} - sG_s(y_{k-1}),                                                                                                                     \label{eqn: m-fista1-gradient}       \\
& x_{k} = \left\{ \begin{aligned} 
                          & z_{k-1}, \quad \text{if}\; f(z_{k-1}) \leq f(x_{k-1}),  \\
                          & x_{k-1}, \quad \text{if}\; f(z_{k-1}) > f(x_{k-1}),
                          \end{aligned} \right.                                                                                                                    \label{eqn: m-fista1-comparison}   \\    
& y_{k} = x_{k} + \frac{(k-1)^{\alpha}}{k^{\alpha}+rk^{\alpha-1}} \left( x_{k} - x_{k-1} \right) + \frac{(k-1)^{\alpha} + r(k-1)^{\alpha-1}}{k^{\alpha}+rk^{\alpha-1}} \left(z_{k-1} - x_{k} \right),            \label{eqn: m-fista1-momentum}
\end{empheq}
\end{subequations}
where the s-proximal subgradient $G_s(y_k)$ replaces the gradient $\nabla f(y_k)$ in \texttt{NAG}-$\alpha$. This formulation explicitly accounts for the non-smooth nature of the composite objective function $\Phi$, while maintaining a similar structure to \texttt{M-NAG}-$\alpha$.

To establish the controllable $O\left( 1/k^{2\alpha} \right)$ convergence rate for \texttt{M-FISTA}-$\alpha$, we generalize the Lyapunov function~\eqref{eqn: mnag-lyapunov} used in~\texttt{M-NAG}-$\alpha$ to accommodate the composite objective function $\Phi$. The new Lyapunov function takes the following form:
\begin{align}
\mathcal{E}(k) = & sk^{\alpha} \left( k^{\alpha} + rk^{\alpha-1} \right)\left( \Phi(x_{k}) - \Phi(x^{\star}) \right)   \nonumber \\
                           &+ \frac12\left\| \left[ k^{\alpha} + r\left( k^{\alpha-1} - (k-1)^{\alpha-1} \right) \right] \left( y_{k} - x_{k}\right)  +  r(k-1)^{\alpha - 1}\left( y_k - x^{\star} \right) \right\|^2. \label{eqn: mfista-lyapunov}
\end{align}
In~\Cref{subsec: m-nag}, the derivation of the $O\left( 1/k^{2\alpha} \right)$ rate for smooth functions hinges on two key inequalities,~\eqref{eqn: fund-inq-xk-xk-mnag} and~\eqref{eqn: fund-inq-xk-xstar-mnag}. These inequalities are generalized in the composite setting to accommodate the proximal components, leveraging~\Cref{lem: fund-inq-composite}, which extends the fundamental inequality for strongly convex functions to the proximal case. These adaptations allow for a rigorous analysis of~\texttt{M-FISTA}-$\alpha$, effectively addressing the non-smooth components of $\Phi$. The convergence properties of~\texttt{M-FISTA}-$\alpha$ are formally encapsulated in the following theorem:
\begin{theorem}
\label{thm: m-fista-alpha}
\label{thm: fista-alpha}
Let $f \in \mathcal{S}_{\mu, L}^1(\mathbb{R}^d)$ and $g \in \mathcal{F}^{0}(\mathbb{R}^d)$, and suppose the step size satisfies $s \in (0, 1/L]$. If the controllable parameter is chosen such that $r > 2\alpha$, then there exists a positive integer $K:= K(\alpha, r)$ such that the iterative sequence $\{x_k\}_{k=0}^{\infty}$, generated by \texttt{FISTA}-$\alpha$ under updates,~\eqref{eqn: m-fista-gradient} ---~\eqref{eqn: m-fista-momentum} with any initial $x_0 = y_0 \in \mathbb{R}^d$, converges at the following rate:
\begin{equation}
\label{eqn: m-fista-alpha-rate}
\Phi(x_k) - \Phi(x^{\star}) \leq \frac{\mathcal{E}(K)}{sk^\alpha \left( k^\alpha + rk^{\alpha-1} \right)}, \quad \forall k \geq K.
\end{equation}
Furthermore, if the step size is set critically as $s = 1/L$, the iterative sequence $\{x_k\}_{k=0}^{\infty}$ satisfies:
\begin{equation}
\label{eqn: m-fista-alpha-rate-l}
\Phi(x_k) - \Phi(x^{\star}) \leq \frac{L\mathcal{E}(K)}{k^\alpha \left( k^\alpha + rk^{\alpha-1} \right)}, \quad \forall k \geq K.
\end{equation}
\end{theorem}

\section{Conclusion and future work}
\label{sec: conclusion}

In this paper, we introduce a family of controllable momentum coefficients for forward-backward accelerated methods, with a focus on the critical step size $s=1/L$. Unlike traditional linear forms, the proposed momentum coefficients adopt an $\alpha$-th power structure, with the parameter $r$ adaptively tuned to $\alpha$. A key contribution is the development of a Lyapunov function tailored to the parameter $\alpha$, which excludes the kinetic energy term. This simplification enabled us to establish a controllable $O\left(1/k^{2\alpha} \right)$ convergence rate for the \texttt{NAG}-$\alpha$ method, provided $r > 2\alpha$. At the critical step size, \texttt{NAG}-$\alpha$ achieves an inverse polynomial convergence rate of arbitrary degree by adjusting $r$ according to $\alpha > 0$. By omitting the kinetic energy term, we further simplified its expression in terms of the iterative sequences $x_k$ and $y_k$, thereby eliminating the need for phase-space representations. This insight allowed us to extend the controllable $O \left(1/k^{2\alpha} \right)$ rate to the monotonic variant, \texttt{M-NAG}-$\alpha$, further enhancing optimization efficiency. Finally, by leveraging the fundamental inequality for composite functions derived in~\citep{li2024linear2}, we extended the controllable $O\left(1/k^{2\alpha} \right)$ rate to proximal algorithms, including  \texttt{FISTA}-$\alpha$ and \texttt{M-FISTA}-$\alpha$.

A notable reference point for first-order methods is the lower bound established in~\citep{nemirovsky1983problem, nesterov2018lectures}, which applies to a finite number of iterations ($k \leq d/2$). In the asymptotic sense, however, the convergence rate achieves $o(1/k^2)$, as shown in~\citep{attouch2016rate}. From our analysis in~\Cref{sec: nag}, particularly~\eqref{eqn: iv-12} and~\eqref{eqn: iv-22}, we observe that the strongly convex conditions used to bound the iterative difference is significantly larger than the corresponding terms,  exceeding them by at least two orders ($O(k^2)$). If we consider a convex objective function augmented with a penalty term that asymptotically approaches zero, such as:
\[
\min_{x\in \mathbb{R}^x} g_k(x)= f(x) + \frac{1}{2k}\|x - x^{\star}\|^2,
\]
then by applying~\texttt{NAG}-$\alpha$ with $r > 2\alpha$, it is possible to achieve the $O\left( 1/k^{2\alpha}\right)$ convergence rate for $\alpha \in (1, +\infty)$. However, in practical scenarios, constructing such a form is infeasible, as it presupposes prior knowledge of the optimal solution, which contradicts the essence of iterative optimization. This reveals a significant gap between theoretical results and practical applicability, as the strongly convex conditions used are overly restrictive and rarely hold in real-world problems. Therefore, exploring weaker conditions beyond strong convexity is a promising direction for future research. Such efforts could help bridge the gap between theory and practice, expanding the practical applicability of accelerated optimization methods.

\section*{Acknowledgements}
We thanks Bowen Li for helpful discussions. Mingwei Fu was partially supported by the Hua Loo-Keng scholarship from CAS.  This project was partially supported by the startup fund from SIMIS.

\bibliographystyle{abbrvnat}
\bibliography{sigproc}
\end{document}